\documentclass{article}

\usepackage{amssymb}
\usepackage{amsmath}
\usepackage[dvips]{color}
\usepackage{bbm}
\usepackage{graphics}
\usepackage{graphicx}
\usepackage{pifont}
\usepackage{subfigure}
\usepackage[hmargin=1.4in,vmargin=1.3in]{geometry}

\newcommand{\dd}{\mathrm{d}}
\newcommand{\dx}{\mathrm{d}x}

\newcommand{\al}{\alpha}
\newcommand{\be}{\beta}

\newcommand{\tbe}{\widetilde{\beta}}

\newcommand{\vfi}{\varphi}
\newcommand{\vphi}{\varphi}

\newcommand{\Nfi}{{\cal{N}}_{\vfi}}
\newcommand{\NFi}{{\cal{N}}_{\Phi}}
\newcommand{\Nfie}{{\cal{N}}_{\vfi}^E}
\newcommand{\omfi}{\omega_{\vfi,_{-}}}
\newcommand{\omFi}{\omega_{\Phi,_{-}}}
\newcommand{\omfiplusz}{\omega_{\vfi,_{+}}}
\newcommand{\omFiplusz}{\omega_{\Phi,_{+}}}
\newcommand{\omfinul}{\omega_{\vfi,_{0}}}
\newcommand{\omFinul}{\omega_{\Phi,_{0}}}

\newtheorem{thm}{Theorem}[section]
\newtheorem{lemma}[thm]{Lemma}

\newtheorem{rem}[thm]{Remark}

\newtheorem{exa}[thm]{Example}

\begin{document}

\markboth{Lajos L\'oczi}{Journal of Difference Equations and Applications}


\title{Discretizing the transcritical and pitchfork bifurcations---conjugacy results}


\author{Lajos L\'oczi\\lloczi@cs.elte.hu\\
Department of Numerical Analysis\\ E\"otv\"os Lor\'and University\\ 
P\'azm\'any P. s\'et\'any 1/C, Budapest, H-1117, Hungary\begin{footnote}
{The project was supported by the European Union and co-financed by the
European Social Fund (grant agreement no.\ TAMOP
4.2.1/B-09/1/KMR-2010-0003). Support from Award No.\ FIC/2010/05 -- 2000000231 made by King Abdullah University of Science and Technology (KAUST) is also acknowledged.}
\end{footnote}\footnote{The final version of this manuscript is to appear in \textit{Journal of Difference Equations and Applications}}}

\maketitle

\begin{abstract}
We present two case studies in one-dimensional dynamics concerning the discretization of transcritical (TC) and pitchfork (PF) bifurcations. In the vicinity of a TC or PF bifurcation point and
under some natural assumptions on the one-step discretization method of order $p\ge 1$, we show that the time-$h$ exact and the step-size-$h$ discretized dynamics are topologically equivalent
by constructing a two-parameter family of conjugacies in each case. 
As a main result, we prove that 
the constructed conjugacy maps are ${\cal{O}}(h^p)$-close to the identity and these estimates
are optimal.


\end{abstract}

\section{Introduction}\label{introduction_tk}

Let us consider a one-dimensional ordinary differential equation 
\begin{equation}\label{kozdiff_tk}\dot{x}(t)=f(x(t),\alpha)\end{equation}
depending on a  scalar bifurcation parameter $\alpha \in \mathbb{R}$. 
Suppose that a one-step discretization of (\ref{kozdiff_tk}) is given by
\begin{equation}\label{diszkretizalt_tk}
X_{n+1}:=\varphi(h,X_n,\alpha),\quad \quad n\in\mathbb{N}\equiv\{0,1,2,\ldots\},
\end{equation}
where $h>0$ denotes the step-size of the 
method $\varphi:\mathbb{R}^+\times\mathbb{R}\times \mathbb{R}\rightarrow \mathbb{R}$
of order $p \in\mathbb{N}^+$, and sufficient smoothness is assumed on both $f:\mathbb{R}\times \mathbb{R}\rightarrow \mathbb{R}$ and $\varphi$. We suppose that (\ref{kozdiff_tk}) undergoes either a transcritical or a pitchfork bifurcation at parameter value $\alpha=0$ with the bifurcating equilibrium located at $x=0$. Let $\Phi(h,\cdot,\alpha):\mathbb{R}\rightarrow \mathbb{R}$ denote the time-$h$-map of the exact solution
flow induced by (\ref{kozdiff_tk}) at parameter value $\alpha$. Then, by
the definition of the order of the numerical method, we have  
\begin{equation}\label{kozelseg_tk}
\left|\Phi(h,x,\alpha)-\varphi(h,x,\alpha)\right|\le const\cdot h^{p+1},\quad \forall\, h\in [0,h_0], 
\forall\, |x|\le \varepsilon_0, \forall\, |\alpha|\le \alpha_0
\end{equation}
with some positive constants $const$, $h_0$, $\varepsilon_0$ and $\alpha_0$, depending on $f$ and $\varphi$.

A fundamental task in numerical analysis is to compare the exact dynamics with the discretized one,
that is, the discrete dynamical systems generated by the iterates of $\Phi$ and $\varphi$, respectively,
both from a qualitative and a quantitative point of view. Let us recall some specific previous results in this direction which served as a direct motivation for our present work. 

The case when $\dot{x}(t)=f(x(t))$ has a hyperbolic equilibrium has been solved satisfactorily. 
It has been proved that one-step methods faithfully 
reproduce the phase portrait near hyperbolic equilibria in the following sense: 
the time-$h$-map of the exact dynamics and the step-size-$h$-map of the 
discretized dynamics are conjugate---i.e., there exists a $C^0$ (and in general) non-linear coordinate transformation
between them---moreover, the conjugacy map can be chosen to be ${\cal{O}}(h^p)$-close to the identity. This type of result can be recast by saying that the original system is numerically structurally stable.

A natural next step is to investigate the relationship between exact and discretized dynamics when hyperbolicity is violated, for example, in the context of bifurcating one-parameter families of ODEs, like (\ref{kozdiff_tk}). Simple examples show that the original dynamics and its discretization need not be conjugate near a general non-hyperbolic equilibrium. Nevertheless, by extending some earlier results, we showed in \cite{garayloczi} that the time-$h$-map of the exact dynamics and the step-size-$h$-map of the 
discretized dynamics are again conjugate in a neighborhood of a fold bifurcation point. In addition, we proved optimal ${\cal{O}}(h^p)$-closeness results in the $\alpha\le 0$ region where the two branches of fixed points are located and merge at $(\alpha,x)=(0,0)$; 
in the fixed-point-free $\alpha>0$ region however we only managed to get
either a singular ${\cal{O}}(h^p \, |\ln \alpha |)$  or a ${\cal{O}}(h)$ estimate between the conjugacy map $J$ and the identity $id$, depending on the construction of $J$. A detailed proof of the $\alpha\le 0$ case is given in \cite{garayloczi}, while a proof of the singular estimate is published in \cite{paezloczisurvey}. Therefore, although numerical structural stability has been established, the problem of optimal conjugacy estimates in the fold bifurcation case is still open.

\subsection{An outlook to the relevant literature }

There is a rich literature on connections between discretizations and bifurcations. 
The interested reader will find a selection of about 50 references on this topic in our
recent work \cite{paezloczisurvey}.
Let us add that the above brief summary in the previous paragraphs is also put into a wider context 
 in \cite{paezloczisurvey}, including the notion of structural and numerical structural stability with earlier important results under various hyperbolicity conditions; references on the properties of numerical methods near bifurcation points of ODEs together with convergence questions; earlier proofs on the discretized fold bifurcation; finally, several theorems on other types of discretized bifurcations of codimension 1, 2 or higher, in 1 or $N$ dimensions, and comparisons between the corresponding bifurcation diagrams, critical eigenvectors or normal form coefficients. These last results are more general, but weaker than conjugacy results. To the best of our knowledge, conjugacy results
in a non-hyperbolic setting have only been obtained for the fold bifurcation.

\subsection{The main results of the present paper}

The purpose of the present work is to establish conjugacy results for
the TC and PF bifurcations. These results are based on the author's thesis \cite{dissertation}.
We show that in a neighborhood of a TC or PF point the original system (\ref{kozdiff_tk}) is not numerically structurally stable if arbitrary one-step methods are considered, but when the allowed $\varphi$ discretizations are suitably restricted to a certain class (containing all Runge--Kutta methods, for example), then numerical structural stability is recovered near TC and PF points, too. 
More specifically, we prove in the TC and PF cases that for any fixed $h$ and $\alpha\in\mathbb{R}$ (with
$h>0$ and $|\alpha|$ sufficiently small) and under some natural assumptions on $\varphi$, there is a solution
$J(h,\cdot,\alpha)$ to the conjugacy equation
\begin{equation}\label{konjugacios alapegyenlet}
J(h,\Phi(h, x,\alpha),\alpha)=\varphi(h,J(h,x,\alpha),\widetilde{\alpha}),
\end{equation}
with $\alpha \mapsto {\widetilde \alpha}$ $(={\widetilde \alpha}(h))$ being a homeomorphism and
 $|\alpha -  {\widetilde
\alpha}|={\mathcal{O}}(h^p)$.\footnote{Throughout the article, the 
constants in the ${\mathcal{O}}(h^p)$ symbols are independent of
$x$ and $\alpha$.}  
We remark that aligning the bifurcation parameter $\alpha$ with $\widetilde{\al}$ in this functional  equation (\ref{konjugacios alapegyenlet}) is usually necessary, since numerical methods may shift the original bifurcation point $x=0\in\mathbb{R}$.
It is also important that the constructed homeomorphism
$J(h,\cdot,\alpha)$ is defined on a fixed, non-shrinking neighborhood
(being independent of the parameters as $h\to 0^+$ or $\alpha\to 0$) of the 
point $x = 0$.  
Beyond the above results we prove in both cases that the distance between the identity map and
the  constructed conjugacy satisfies an  
\[|id - J(h,\cdot, \alpha)|={\cal{O}}(h^p)\] closeness estimate, and this estimate is optimal in $h$.

\subsection{The organization of the paper}

The structure of the paper is as follows.  Section \ref{notationsection} summarizes our notation. 
The TC bifurcation is treated in Section \ref{TCmainresultssection} (giving an overview
of the results) and in Section \ref{conjugacyTC} (containing the detailed proofs).
In Section \ref{TCmainresultssection},  the defining conditions for the TC bifurcation are first 
analyzed together with our assumptions on the discretization $\varphi$. Necessity of these assumptions is illustrated by an example. The normal forms of $\Phi$ and $\varphi$ and their closeness 
relations are discussed next. With the help of these normal forms and 
the corresponding fundamental domains, the construction of the conjugacy map is 
sketched.  Section \ref{TCmainresultssection} then concludes with our main theorem
on the TC conjugacy estimates. The detailed construction of the conjugacy is 
actually presented in  Section \ref{tk_constructionoftheconjugacy}, proving topological equivalence between the original time-$h$-map and its discretized counterpart. The quantitative Sections \ref{tk_optimality}--\ref{tk_outerregion} deal with the optimality of the ${\mathcal{O}}(h^p)$ estimates
first, then with the explicit conjugacy estimates in the inner and outer regions in the $(\alpha,x)$-plane, respectively. As for the PF bifurcation, the arrangement of the material is analogous. Section \ref{PFmainresultssection} gives a summary of
our results in the PF case, Section \ref{pf_constructionoftheconjugacy} describes the construction 
of the conjugacy in detail, while Sections \ref{pf_optimality}--\ref{pf_outerregion} 
contain the optimal ${\mathcal{O}}(h^p)$ PF closeness estimates. 
Finally in Section \ref{Section5}, we provide a numerical example demonstrating the closeness
estimates in the TC case, and mention some open questions about a possible generalization of
the above results to higher dimensions.

\subsection{Notation}\label{notationsection}

The floor and ceiling functions are $\lfloor \cdot \rfloor$ and 
$\lceil \cdot \rceil$, as usual. For $a$, $b \in \mathbb{R}$, $[\{a,b\}]$ is the closed interval generated by the set $\{a,b\}$, that is  $[\{a,b\}]:=[\min(a,b),\max(a,b)]$. 

The $k^\mathrm{th}$ iterate ($k\in\mathbb{Z}$) of a real function $g$ is denoted by $g^{[k]}$. When $k$ is negative, $g$ is assumed to be invertible, so, for example, $g^{[-1]}$ denotes the inverse 
function of $g$.

By ``smooth function" we mean a $C^k$ function (in all variables) with sufficiently large $k\in\mathbb{N}^+$
 and with the last derivative bounded.

The normal form of $\Phi(h,\cdot,\al)$ is denoted by $\NFi(h,\cdot,\al)$, while $\Nfi(h,\cdot,\al)$ denotes the normal form of $\varphi(h,\cdot,\al)$. Iterates of one of the normal forms with suitable starting values will define the real sequences $x_n$, $y_n$ or $z_n$, used in the definition of the fundamental domains (their dependence on $h$ and $\al$ is often suppressed however). 

The symbols $\omFi$, $\omFinul$ and $\omFiplusz$ stand for certain fixed points of the maps
$\Phi(h,\cdot,\al)$ (or, of the $\NFi(h,\cdot,\al)$ normal forms). These fixed points may depend on $h$ and $\al$. The meaning of $\omfi$, $\omfinul$ and $\omfiplusz$ is analogous.

As already mentioned earlier, $h_0>0$, $\varepsilon_0>0$ and $\alpha_0>0$ are (small) positive constants, depending on $f$ and $\varphi$. Some restrictions on their values will be given later. The particular positive constant in the estimate of the distance of the normal forms is denoted by $c$
(see (\ref{tk_normalformakkozelsege}) and (\ref{pf_normalformakkozelsege})). It is independent of $h$, $x$ and $\al$, and its value is fixed within Section \ref{conjugacyTC} or \ref{conjugacyPF}. A generic positive constant (independent of $h$, $x$ and $\al$) is denoted by $const$. It may denote different positive numbers at different occurrences. We will use $K$ to denote a positive uniform bound on the moduli of the functions $\eta_3$, $\widetilde{\eta}_3$ or
$\eta_4$, $\widetilde{\eta}_4$ together with their first and second derivatives. These $\eta$ functions will appear in the tails of the normal forms. The value of $K$  is assumed to be fixed within Section \ref{conjugacyTC} or \ref{conjugacyPF}.

Various (mixed) partial derivatives with respect to the corresponding variables are denoted by 
$f_x$, $f_{x\al}$, and so on. We use superscript $^E$ for function evaluation at general (but fixed) parameter values $h$ and $\al$, so $J^E$, for example, abbreviates the function $J(h,\cdot,\al)$, while evaluation at the bifurcation point is denoted by the $^B$ operator: $f_x^B$, for example, stands for $f_x(0,0)$.

Finally, let us call the reader's attention to a slight notational incompatibility  in bifurcation theory regarding the order of variables: bifurcation diagrams are usually depicted in the $(\al,x)$-plane,
however, as arguments of a function, the variable $x$ generally precedes the
bifurcation parameter $\al$.

\section{Main Results}

\subsection{The TC case}\label{TCmainresultssection}

Suppose that the origin $(x,\alpha)=(0,0)\in\mathbb{R}^2$ is an equilibrium as well as a transcritical bifurcation point for
(\ref{kozdiff_tk}), that is, the following conditions hold:
\begin{subequations}\label{transcritical conditions}
\begin{align} 
        f(0,\al)=0,\quad \forall\, |\alpha|\le \alpha_0, \\
        f_{x}^B=0,\quad f_{xx}^B\ne 0,\quad f_{x\alpha}^B\ne 0.
\end{align}
\end{subequations}

\begin{rem}\label{tkintroductoryremark}
 Besides (\ref{transcritical conditions}), there are other ways to formulate the TC conditions in one dimension for differential equations. Similarly, various 
sets of conditions can be given for maps guaranteeing the appearance of a TC
bifurcation. 
For a map $x\mapsto g(x,\al)$ to undergo a TC bifurcation near the origin, it is sufficient to have
\begin{equation}\label{TCsuffconditionsinterval}
g(0,\al)=0,\quad \forall\, |\alpha|\le \alpha_0,
\end{equation}
\begin{equation}\label{TCsuffconditionspoint}
g_{x}^B=1,\quad g_{xx}^B\ne 0, \quad g_{x\alpha}^B\ne 0.
\end{equation}
Notice that unlike conditions (\ref{TCsuffconditionspoint}), (\ref{TCsuffconditionsinterval}) is 
not a point condition on $g$ at the bifurcation point $(x,\alpha)=(0,0)$, so 
one can ask whether (\ref{TCsuffconditionsinterval}) can be relaxed. The answer is affirmative,
but some care should be exercised. 
In  \cite{wiggins}, for example, instead of 
(\ref{TCsuffconditionsinterval})--(\ref{TCsuffconditionspoint}), the map $g$ is simply required to satisfy 
\begin{equation}\label{wiggdefcond}
g^B=0,\quad g_\al^B=0,\quad g_{x}^B=1,\quad g_{xx}^B\ne 0,\quad g_{x\alpha}^B\ne 0.
\end{equation}
The argument in \cite{wiggins} 
tacitly assumes in addition that $x$ can be factored out from $g$,
and prove from these the presence of a TC-structure near the origin.
However, the example $x_{n+1}:=g(x_n,\al)$ with
$g(x,\al):=\al^2+(1+\al)x+x^2$ illustrates that (\ref{wiggdefcond}) alone 
is insufficient to define a TC bifurcation: 
since $(x,\al)=(0,0)$ is the only fixed point of this map, clearly no 
TC bifurcation of fixed points can occur here.  
It is unfortunate that some other mathematical works, e.g., \cite{fosterkhumalo}, or teaching materials also try to define the TC bifurcation imperfectly as (\ref{wiggdefcond}).
Nevertheless, TC bifurcation for maps can be guaranteed via point conditions only. In \cite{glendinning}, for example, a discriminant condition is used: conditions
\[
g^B=0,\quad g_\al^B=0,\quad g_{x}^B=1,\quad g_{xx}^B\ne 0,\quad \left(g_{x\alpha}^B\right)^2-g_{xx}^B\cdot g_{\al\al}^B>0
\]
imply the presence of a TC bifurcation near the origin. 
To complicate the matter a bit, \cite{glendinning} contains a typo (but the intended meaning is
clear from the context):
instead of the last $>$ inequality above, a ``$\ne$" sign appears there. To summarize, conditions
(\ref{transcritical conditions}), or (\ref{TCsuffconditionsinterval})--(\ref{TCsuffconditionspoint}) we adopted are not the weakest ones, but they are simple and still retain all essential features of the problem. Finally, it is instructive to compare this remark with Remark \ref{pfintroductoryremark}.
\quad $\Box$
\end{rem}

To ensure that the origin $x=0$ is a fixed point also for the discretization map
\begin{equation}\label{mapphi_tk}x\mapsto\varphi(h,x,\al),\end{equation}  
we assume that 
\begin{equation}\label{phiequilibrium_tk}
\vphi(h,0,\al)=0
\end{equation} 
holds for all sufficiently small $h\ge 0$ and $|\al|$.  This condition is necessary for (\ref{mapphi_tk}) to undergo a TC bifurcation near the origin, as illustrated by the following example.\\

\begin{exa} Let us fix some $h>0$. Suppose we have a map $\vphi(h,x,\al):=h^{2p+1}+(1+h\al)x+h x^2$,
violating (\ref{phiequilibrium_tk}), but satisfying (\ref{kozelseg_tk}) with  
$\Phi(h,x,\al):=(1+h\al)x+h x^2$. This $\Phi$ undergoes a TC bifurcation near the origin, however, $\vphi$ does not: its fixed points are given by $x_{\pm}=\frac{1}{2} \left(-\alpha \pm \sqrt{\alpha ^2-4h^{2p}}\right)$, so the fixed-point branches of $\vphi$ are separeted by the
central $\al\in (-2h^p,2h^p)$ strip.\quad $\Box$
\end{exa}

\begin{rem}\label{remark2} It is well-known that all Runge--Kutta methods preserve equilibria (see
Lemma \ref{pfRKlemmainonedimension}), hence (\ref{phiequilibrium_tk}) is automatically satisfied for these discretizations.\quad $\Box$
\end{rem}

Now let us turn our attention to the solution of equation (\ref{konjugacios alapegyenlet}). 
Since it will be easier to construct $J$ and prove the closeness estimates 
by working
with some concrete normal forms instead of $\Phi$ and $\varphi$, we present 
two results about the structure of the normal forms for the maps $\Phi(h,\cdot,\al)$ and 
$\varphi(h,\cdot,\al)$  
near the equilibrium being also a TC bifurcation point.

\begin{lemma}\label{Phinormalforma_tk} Suppose that $f\in C^{p+6}$ and (\ref{transcritical conditions}) hold. Then there are smooth and invertible coordinate and parameter changes transforming the map
$x\mapsto\Phi(h,x,\al)$
into
\[\eta\mapsto (1+h\beta)\eta+s\cdot h\eta^2+h\eta^3\cdot \eta_3(h,\eta,\beta),\]
where $\eta_3$ is a 
smooth function, and $s$ is either $1$ or $-1$ for all small values of $h\ge 0$ and $|\alpha|$.
\end{lemma}

\begin{thm}\label{phinormalforma and closeness_tk} Suppose that $f\in C^{p+6}$ and conditions (\ref{kozelseg_tk}), (\ref{transcritical conditions}) and (\ref{phiequilibrium_tk}) hold with
a sufficiently smooth $\varphi$. Then there are smooth and invertible coordinate and parameter changes transforming the map
$x\mapsto\vphi(h,x,\al)$
into
\[\widetilde{\eta}\mapsto (1+h\widetilde{\beta})\widetilde{\eta}+s\cdot 
h{\widetilde{\eta}}^2+h\widetilde{\eta}^3\cdot
\widetilde{\eta}_3(h,\widetilde{\eta},\tbe),\] 
where $\widetilde{\eta}_3$ is a smooth function, and $s$ is either $1$ or $-1$ for all small values of $h\ge 0$ and $|\alpha|$. Moreover, the smooth invertible coordinate and parameter changes above and those in Lemma \ref{Phinormalforma_tk} are $\cal{O}$$(h^p)$-close to each other, and 
\[
|\eta_3-\widetilde{\eta}_3|\le const\cdot h^p.
\]
\end{thm}

We omit the technical proofs of the above two results; they are given in \cite{dissertation}. The  idea of their proofs is to follow the corresponding sections of
\cite{kuznetsov} with $h$ suitably built into the computations:
the smoothness assumptions allow multivariate Taylor polynomials of
$\Phi$ and $\varphi$ with integral remainders to be constructed
explicitly, then qualitative and quantitative inverse and implicit
function theorems (see, e.g.,  \cite[Theorem 2.2]{paezloczisurvey}) are applied.
We remark that the value of $s=\pm 1$ is the same in Lemma \ref{Phinormalforma_tk} and Theorem
\ref{phinormalforma and closeness_tk}, and this sign is only responsible for the orientation of the bifurcation diagrams. For definiteness, from
now on in this section and in Section \ref{conjugacyTC}, $s:=1$
is chosen.

Finally, \cite{dissertation} also proves that 
$|\be-\tbe|\equiv |\be(h,\al)-\tbe(h,\al)|\le const\cdot h^{p}$, so
we can apply a further parameter shift $\tbe \mapsto \be$ in the normal form in Theorem \ref{phinormalforma and closeness_tk}, being $\cal{O}$$(h^p)$-close to the identity. 
This allows us to use the    
bifurcation parameter $\al$ again instead of $\be$ and $\tbe$. 
Moreover, as a notational simplification, we will write $x$ instead of 
the dummy variables $\eta$ and $\widetilde{\eta}$.\\

The above reductions imply that it is enough to
solve the conjugacy equation (\ref{konjugacios alapegyenlet}) for
the corresponding normal forms and with $\alpha =  {\widetilde \alpha}$. 
This will be carried out in Section \ref{tk_constructionoftheconjugacy} by
using the method of fundamental
domains, see \cite{kuznetsov} also. The essence of this method is the following:
for fixed $h$ and $\alpha$, $J(h, \cdot, \alpha)$ is
prescribed on a suitable starting interval and then recursively
extended by using rearrangements of (\ref{konjugacios alapegyenlet})
with $\alpha =  {\widetilde \alpha}$. Near attractive fixed points,
the fundamental domains are defined by forward iterates
\begin{equation}\label{fund_dom_seq}
x_{k+1}(h,\alpha):={\cal{N}_{\vfi}}(h,x_k,\al) \quad (k=0,1,2,\ldots)
\end{equation}
of a suitable starting value $x_0$,
while backward iterates are used near repelling fixed points.
Some members of the sequences $x_k(h,\alpha)$ are illustrated in Figure \ref{figure3}. 
From the proposed construction in Section \ref{tk_constructionoftheconjugacy}, it is easily seen that
$J(h,\cdot,\alpha)$ is a homeomorphism locally near the origin,
moreover, $J$ is continuous in its first and third arguments as well.\\ 

After the conjugacy map has been constructed, estimating 
 $|x-J(h,x,\alpha)|$ from above for $0<h\le h_0$,
$|x|\le\varepsilon_0$ and $|\alpha|\le \alpha_0$ is much more
involved. This is described in Sections \ref{tk_optimality}--\ref{tk_outerregion}.

Section \ref{tk_optimality} illustrates by simple examples that the distance of the 
fixed points of
$\NFi(h,\cdot,\al)$ and $\Nfi(h,\cdot,\al)$  can be bounded 
from below by $\mathit{const}(\alpha)\cdot h^p$. Since a
conjugacy necessarily maps fixed points into fixed points, better
estimates than ${\cal{O}}(h^p)$ for $|x-J(h,x,\alpha)|$ generally
cannot be expected---at least in terms of $h$.

Then, by using the
recursive definitions of the conjugacy, it turns out that for an
${\cal{O}}(h^p)$ closeness estimate, a discrete Gronwall-type estimate
of the form
\begin{equation}\label{introduction_kozelsegi_becslesek}
\left(h\sum_{k=0}^{n}|x_k|^\omega
\left(\prod_{j=k}^{n} \left({\cal{N}}_{\Phi}\right)_x(h,x_k,\alpha)\right)\right)\cdot
h^p\le \mathit{const}\cdot h^p
\end{equation}
suffices, with $\omega=3$, the sequence $x_k$ defined by
(\ref{fund_dom_seq}), and $\mathit{const}$ being independent of
$0<h\le h_0$, $|x|\le\varepsilon_0$ and $|\alpha|\le \alpha_0$. The
difficulty in estimating the left-hand side of (\ref{introduction_kozelsegi_becslesek}) is
that the derivative of the normal form is approximately $1$ and it is 
evaluated at an implicitly defined non-linear sequence $x_k$, so the
contribution of the product is not easily established. At this
point, a special family of non-linear recursions of T.~H\"uls
\cite{thorstenapplmath} proved to be invaluable: for $a>0$ and $0<q\in
\mathbb{N}$ arbitrary parameters and $z_0>0$ sufficiently small
starting value, the closed form solution of
\begin{equation}\label{thulsmodelfunction}
z_{n+1}:=\frac{z_n}{(1+a q z_n^q)^{1/q}} \quad (n=0,1,2,\ldots)
\end{equation} is
given by
\[
z_{n}=\frac{z_0}{(1+n a q z_0^q)^{1/q}}.
\]
This family of model functions coupled with the 
power of \textit{Mathematica} made it possible to 
formulate our key lemmas Lemma \ref{tk_xknagysagrend} and Lemma \ref{tk_zknagysagrend},
providing convergence
speed estimates
for the parametric sequence $x_k$ and its ``outer" counterpart $z_k$.
Then Lemma \ref{tk_inner_1} and Lemma \ref{tk_lemma3.11inner2} settle the closeness estimates in the inner region,
while Section \ref{tk_outerregion} in the outer region.

As a summary of these results concerning the TC case, we get the following theorem.

\begin{thm}\label{Conjugacy estimate in the TC bifurcation}  Suppose that the conditions of Theorem \ref{phinormalforma and closeness_tk} hold.
Then the conjugacy $J$ (constructed in Section \ref{tk_constructionoftheconjugacy}) satisfies closeness estimates that are optimal in $h$ near the TC point: there exists
a positive constant $const>0$ such that for all $0<h\le h_0$, $|x|\le \varepsilon_0$ and $|\alpha|\le \alpha_0$ 
\begin{equation}\label{Theorem2.4closenesspart}
|\,x-J(h,x,\alpha)|\le const\cdot h^p.
\end{equation}
\end{thm}

Let us add finally that it is seen from the proofs that instead of $\textit{const}\cdot h^p$ in the above closeness estimates, sometimes a better $\textit{const}\cdot
 \alpha^\kappa\,h^p$ can be written with $\kappa\in\{1,2\}$. Generally, the closer $x$ lies to the   fixed point branches, the larger $\kappa$ one can choose.

\subsection{The PF case}\label{PFmainresultssection}

Suppose that the origin $(x,\alpha)=(0,0)\in\mathbb{R}^2$ is an equilibrium as well as a pitchfork bifurcation point 
for
(\ref{kozdiff_tk}), that is, the following conditions hold:
\begin{subequations}\label{pitchfork conditions}
\begin{align} 
        f(0,\al)=0,\quad \forall\, |\alpha|\le \alpha_0, \\
	f_{x}^B=0,\quad f_{xx}^B= 0,\quad f_{xxx}^B\ne 0,\quad f_{x\alpha}^B\ne 0.
\end{align}
\end{subequations}

\begin{rem}\label{pfintroductoryremark}
In this context, $f$ is usually assumed to be odd, that is $f(x,\al)=-f(-x,\al)$ for $|x|\le \varepsilon_0$ and $|\alpha|\le \alpha_0$. The above asymmetric PF conditions  (\ref{pitchfork conditions})  are thus weaker than the usual ones. If symmetry were assumed, then the normal form transformations in Lemma \ref{Phinormalforma_pf} and Theorem \ref{phinormalforma and closeness_pf} below would be much easier. Symmetry, however, is not essential here.
We remark that even (\ref{pitchfork conditions}) can be slightly weakened (cf. Remark \ref{tkintroductoryremark}). Regarding the PF conditions for maps, 
 \cite{wiggins}, for example, uses again only point conditions at the bifurcation point $(0,0)$:
unlike in the TC case, this time they are sufficient to define a PF near the origin, as shown by Theorem \ref{pf_asymm} below. Summarizing, similarly to the TC case, we have not adopted the most general conditions for a PF bifurcation in (\ref{pitchfork conditions}), nevertheless symmetry is not assumed.
\quad $\Box$
\end{rem}


\begin{thm}\label{pf_asymm} Assume we have a smooth map $x\mapsto g(x,\al)$  defined near the origin $(0,0)$ and depending smoothly also on $\al$ such that \[g(0,0)=0,\ g_{x}(0,0)=1,\ g_{xx}(0,0)=0,\ g_{xxx}(0,0)\ne 0,\  g_{\alpha}(0,0)=0,\  g_{x\alpha}(0,0)\ne 0,\] 
\[g_{xxx}(0,0)\cdot g_{x\al}(0,0)>0.\] Then the map $g$ undergoes a pitchfork bifurcation locally at the origin: $g$ has precisely 3 branches of fixed points $\rho_0(\al)$ and $\rho_\pm(\al)$
with the following properties. There exist positive constants $c_0>0, c_2>c_1>0$ and $\al_0>0$
such that \[|\rho_0(\al)|\le c_0|\al| \quad\text{ for } |\al|\le \al_0,\] 
\[c_1 |\al|^{1/2}\le \rho_+(\al)\le c_2 |\al|^{1/2} \quad\text{ for } -\al_0\le \al<0,\] and 
\[-c_2 |\al|^{1/2}\le \rho_-(\al)\le -c_1 |\al|^{1/2} \quad\text{ for } -\al_0\le \al<0.\] The case $g_{xxx}(0,0)\cdot g_{x\al}(0,0)<0$ yields the ``mirror-symmetrical" counterpart: there are three branches of fixed points for $\al>0$ and a unique branch for $\al\le 0$ with similar estimates.
\end{thm}

\noindent The elementary proof of the above theorem on asymmetric pitchfork bifurcations can be found in \cite{dissertation} and it is omitted here.\\

As a preparation for the normal form transformations for the maps 
$x\mapsto\Phi(h,x,\al)$ and 
\begin{equation}\label{mapphi_pf}
x\mapsto\varphi(h,x,\al)
\end{equation}  
near the equilibrium being also a PF bifurcation point, let us make the following considerations. To ensure that the origin $x=0$ is a fixed point also for the discretization map (\ref{mapphi_pf}), we assume that 
\begin{equation}\label{phiequilibrium_pf}
\vphi(h,0,\al)=0
\end{equation} 
holds for sufficiently small $h\ge 0$ and $|\al|$. We also suppose that
\begin{equation}\label{phipitchforknecessary_pf}
\vphi_x(h,0,0)=1\quad \textrm{and}\quad \vphi_{xx}(h,0,0)=0
\end{equation} 
hold for all sufficiently small $h\ge 0$. These conditions are necessary and sufficient for (\ref{mapphi_pf}) with property (\ref{kozelseg_tk}) to undergo a PF bifurcation near the origin. Necessity is illustrated by three examples below, while sufficiency is proved by the normal form transformations themselves.

\begin{figure}
\begin{center}
\subfigure{\includegraphics{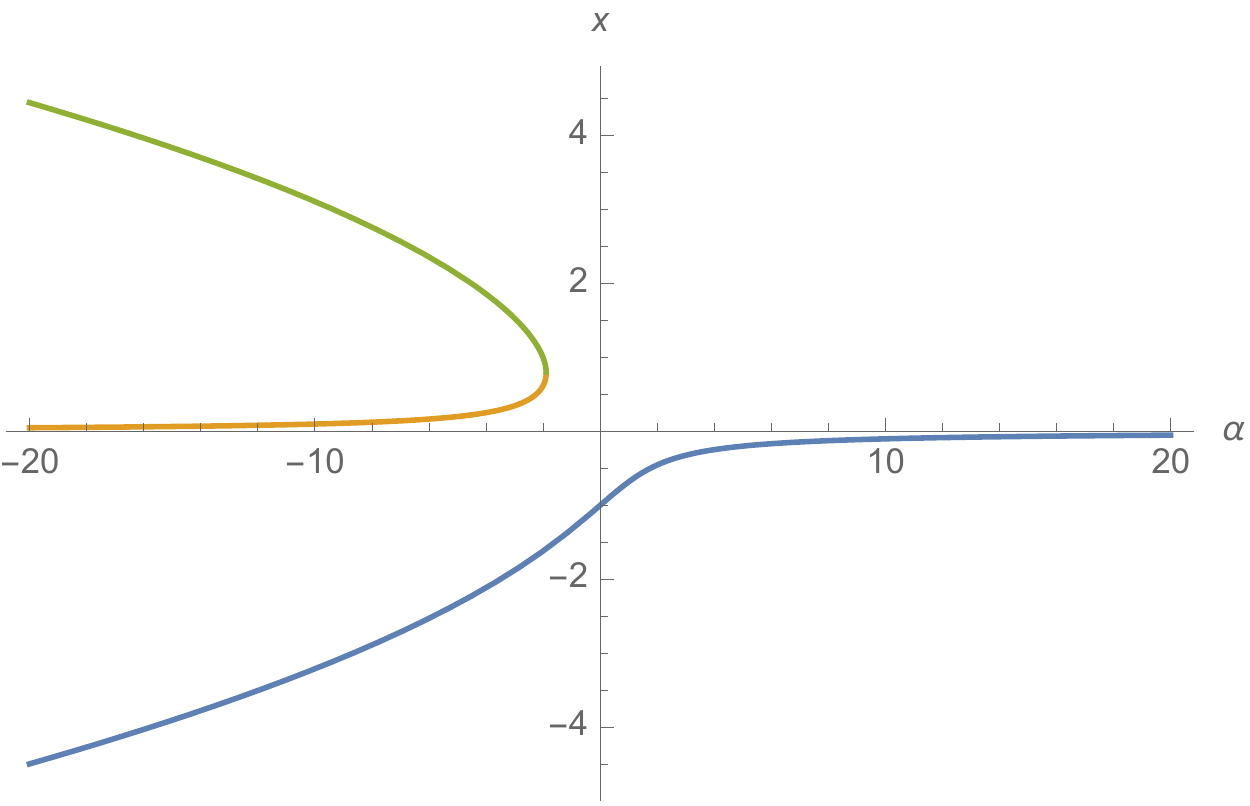}}
\subfigure{\includegraphics{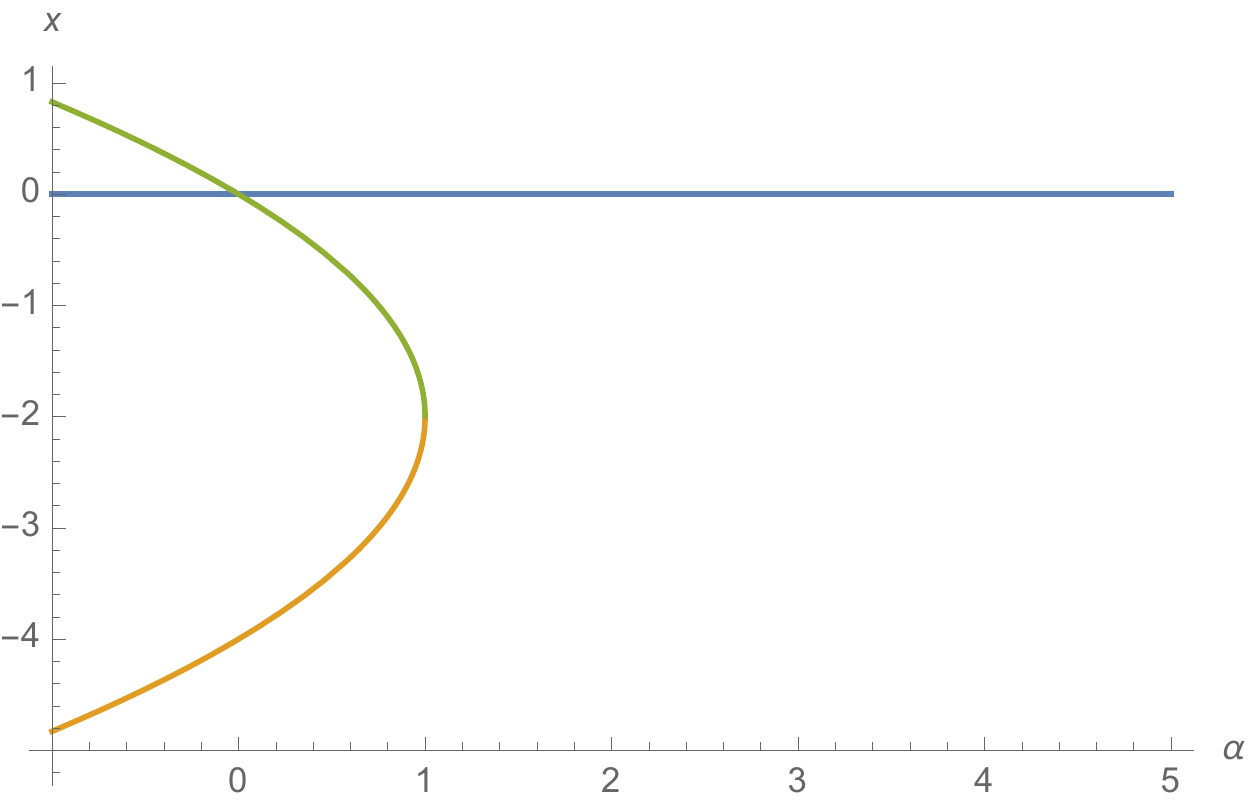}}
\caption{The solution of $\vphi(h,x,\al)=x$ from Example \ref{firstPFexample} 
(upper figure) and Example \ref{secondPFexample} (lower figure) with some  fixed $h>0$.}\label{fig1}
\end{center}
\end{figure}

\begin{exa}\label{firstPFexample} 
Consider the map $\vphi(h,x,\al):=h^{3p+1}+(1+h\al)x+h x^3$. For this $\vphi$, condition $\vphi(h,0,\al)=0$ does not hold, but $\vphi$ satisfies (\ref{kozelseg_tk}) with $\Phi(h,x,\al):=(1+h\al)x+h x^3$. This $\Phi$ has a PF bifurcation at the origin, $\vphi$ however does not, 
see Figure \ref{fig1}.\quad $\Box$
\end{exa}


\begin{exa} \label{secondPFexample} 
Suppose now we have a map $\vphi(h,x,\al):=(1+h\al)x+h^{p+1}x^2+h x^3$. Clearly, $\vphi$ violates condition $\vphi_{xx}(h,0,0)=0$, but satisfies (\ref{kozelseg_tk}) with $\Phi(h,x,\al):=(1+h\al)x+h x^3$. This $\vphi$ does not have a PF near the origin, 
see Figure \ref{fig1}.\quad $\Box$
\end{exa}

\begin{figure}
\begin{center}
\subfigure{\includegraphics{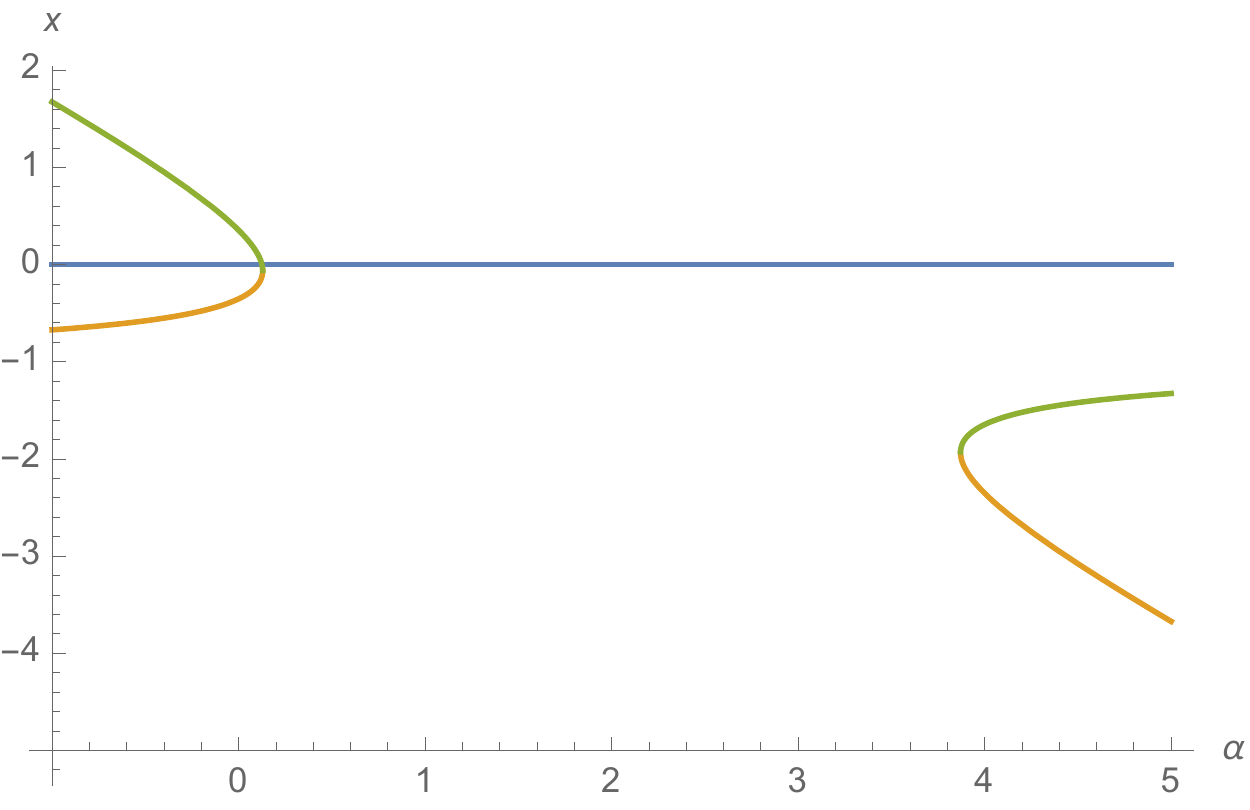}}
\subfigure{\includegraphics{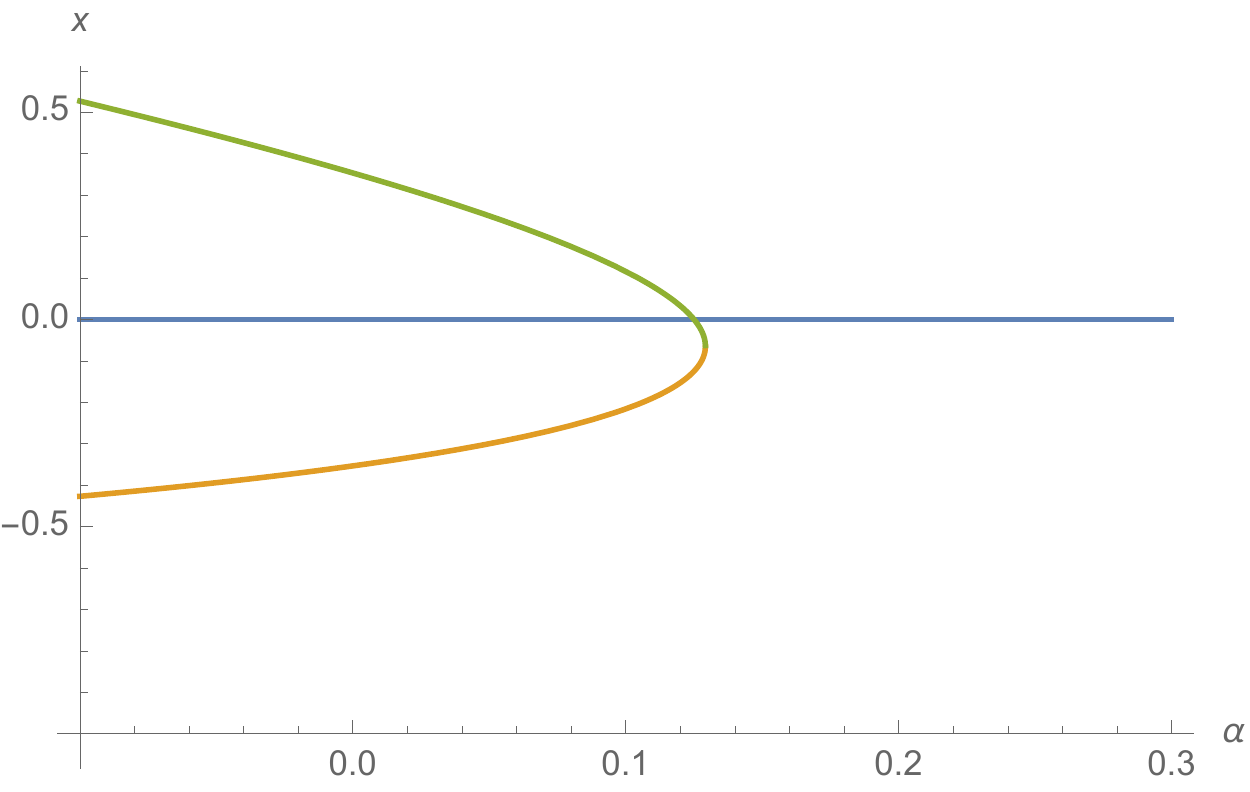}}
\caption{The solution of $\vphi(h,x,\al)=x$ from Example \ref{thirdPFexample} 
(upper figure) with some fixed $h>0$. The lower figure
zooms in on the left half of the upper figure, showing that $\vphi$ does not have a PF bifurcation.}\label{fig2}
\end{center}
\end{figure}

\begin{exa}\label{thirdPFexample} 
Suppose finally that the discretized map has the form 
$\vphi(h,x,\al):=(1+h\al-h^{p+1})x+h \al x^2+h x^3$.
This $\vphi$ does not satisfy $\vphi_x(h,0,0)=1$, but is sufficiently close to
$\Phi(h,x,\al):=(1+h\al)x+h \al x^2+h x^3$. It is easily seen that $\Phi$ has a PF at the origin,
but $\vphi$ does not, see Figure \ref{fig2}.\quad $\Box$
\end{exa}

In the following lemma we prove that all Runge--Kutta methods satisfy the three requirements (\ref{phiequilibrium_pf})--(\ref{phipitchforknecessary_pf}) above. 
\begin{lemma}\label{pfRKlemmainonedimension} Suppose that $f^B=0$, $f_x^B=0$, $f_{xx}^B=0$ and let $\vphi$ be a general $s$-stage Runge--Kutta method with step-size $h$, that is 
\[\vphi(h,x,\alpha) := x+h \sum_{i=1}^s \gamma_{i}\cdot k_i(h,x,\alpha),\]
with each function $k_i$ ($i=1,2,\ldots,s$) satisfying the (implicit) equation 
\[ k_i(h,x,\alpha)=f(x+h\sum_{j=1}^{s} \beta_{i,j}\cdot 
k_j(h,x,\alpha),\alpha),\]
where  $\gamma_i$ and $\beta_{i,j}$ ($i,j=1,2,\ldots,s$) are given real parameters. Then for every $h\ge 0$ sufficiently small, we have that $\vphi_x(h,0,0)=1$ and $\vphi_{xx}(h,0,0)=0$.
\end{lemma}
\textit{Proof.} Since $f^B=0$, from unique solvability we get for all $i$ that $k_i(h,0,0)\equiv 0$, implying the well-known property (\ref{phiequilibrium_pf}). On the other hand, differentiating the implicit defining equation we get that
\[
(k_i)_x (h,0,0)=\left( 1 + h \sum_{j = 1}^{s}\beta_{i,j}\,({k_j})_x(h,0,0) \right)   \,
  f_x\left(h \sum_{j = 1}^{s}\beta_{i,j}\,{k_j}(h,0,0)  ,0\right).
\] 
But $k_j(h,0,0)\equiv 0$ and $f_x(0,0)=0$, hence $(k_i)_x(h,0,0)\equiv 0$, so $\vphi_x(h,0,0)\equiv 1$.

Differentiating the defining equation again we see that
\[
(k_i)_{xx} (h,0,0)=h\sum_{j = 1}^{s}\beta_{i,j}\,({k_j})_{xx}(h,0,0) \,
   f_x\left(h\sum_{j = 1}^{s}\beta_{i,j}\,{k_j}(h,0,0)  ,0\right) + 
\]
\[ \left( 1 + h\sum_{j = 1}^{s}\beta_{i,j}\,({k_j})_x(h,0,0) \right) ^2\,
   f_{xx}\left(h\sum_{j = 1}^{s}\beta_{i,j}\,{k_j}(h,0,0)  ,0 \right),
\]
but $f_x(0,0)=f_{xx}(0,0)=0$, so $(k_i)_{xx} (h,0,0)\equiv 0$ and $\vphi_{xx}(h,0,0)\equiv 0$.
\quad $\Box$

\begin{rem}  It is remarkable that 
$f^B=0$, $f_x^B=0$ and $f_{xx}^B=0$ imply many other vanishing quantities for Runge--Kutta methods. It can be proved in a similar, recursive fashion as above, that, for example,
\[
(k_i)_{h}(h,0,0)\equiv(k_i)_{hh}(h,0,0)\equiv(k_i)_{hx}(h,0,0)\equiv(k_i)_{hxx}(h,0,0)\equiv(k_i)_{hhxx}(h,0,0)\equiv 0.
\]
These proofs require considerably more computations. A consequence, for example, is that \[\vphi_{h x}(h,0,0)\equiv\vphi_{hhxx}(h,0,0)\equiv 0.\] However, these latter two formulae---needed in the proof of Theorem \ref{phinormalforma and closeness_pf} below---directly follow from (\ref{phipitchforknecessary_pf}) as well, so they hold not only for Runge--Kutta methods.\quad $\Box$
\end{rem}

After these introductory remarks, let us present the actual PF normal forms. Since the PF conditions (\ref{pitchfork conditions}) are a special case of the corresponding TC
conditions, we start the normal form transformation just as in the TC case (but truncate Taylor-expansion at fourth order instead of third). 

\begin{lemma}\label{Phinormalforma_pf} Suppose that $f\in C^{p+7}$ and (\ref{pitchfork conditions}) hold. Then there are smooth and invertible coordinate and parameter changes transforming the map
$x\mapsto\Phi(h,x,\al)$ into 
\[
\eta \mapsto (1+h\be)\eta+s\cdot h\eta^3+h\eta^4\cdot \eta_4(h,\eta,\be),\]
where $\eta_4$ is a 
smooth function, moreover, $s$ is either $1$ or $-1$ for all small values of $h\ge 0$ and $|\alpha|$.
\end{lemma}

\begin{thm}\label{phinormalforma and closeness_pf} Suppose that $f\in C^{p+7}$ and conditions  (\ref{kozelseg_tk}), (\ref{pitchfork conditions}), (\ref{phiequilibrium_pf}) and (\ref{phipitchforknecessary_pf}) hold with
a sufficiently smooth $\varphi$. Then there are smooth and invertible coordinate and parameter changes transforming the map
\[x\mapsto\vphi(h,x,\al)\] 
into
\[\widetilde{\eta}\mapsto (1+h\widetilde{\beta})\widetilde{\eta}+s\cdot 
h{\widetilde{\eta}}^3+h\widetilde{\eta}^4\cdot
\widetilde{\eta}_4(h,\widetilde{\eta},\tbe),\] 
where $\widetilde{\eta}_4$ is a smooth function, and $s$ is either $1$ or $-1$ for all small values of $h\ge 0$ and $|\alpha|$. Moreover, the smooth invertible coordinate and parameter changes above and those in Lemma \ref{Phinormalforma_pf} are $\cal{O}$$(h^p)$-close to each other, and
\[
|\eta_4-\widetilde{\eta}_4|\le const\cdot h^p.
\]
\end{thm}

The proofs of the above two results are again given in \cite{dissertation} and we omit them now (we 
just add that obtaining these PF normal forms by using the TC normal forms is still
a non-trivial task). Here again the value of $s=\pm 1$ is the same in the above lemma and theorem, but now and in Section \ref{conjugacyPF}, the value $s:=-1$ is chosen. Then, similarly to the TC case,  we can apply a parameter shift $\tbe \mapsto \be$, being $\cal{O}$$(h^p)$-close to the identity, implying that the bifurcation parameter $\al$ can be used instead of $\be$ and $\tbe$. To simplify the notation, 
we will write $x$ instead of the dummy variables $\eta$ and $\widetilde{\eta}$.

The above reductions again mean that it is enough to
solve the conjugacy equation (\ref{konjugacios alapegyenlet}) for
the corresponding normal forms and with $\alpha =  {\widetilde \alpha}$. 
This is carried out in Section \ref{pf_constructionoftheconjugacy}, again via the
the method of fundamental
domains defined by sequences like (\ref{fund_dom_seq}), of course, now with the corresponding PF normal 
form. It is again true that
$J(h,\cdot,\alpha)$ is a homeomorphism locally near the origin,
and $J$ is continuous in its first and third variables as well.

Sections \ref{pf_optimality}--\ref{pf_outerregion} then describe the conjugacy estimates
in detail. In formula (\ref{introduction_kozelsegi_becslesek}), we have $\omega=4$ in the PF case.
The key results on the convergence speed now are Lemmas \ref{pf_xknagysagrend}
and \ref{pf_zknagysagrend}.

As a conclusion, we obtain the following theorem.

\begin{thm}\label{Conjugacy estimate in the PF bifurcation} Suppose that the conditions of Theorem \ref{phinormalforma and closeness_pf} hold.
Then the conjugacy $J$ (constructed in Section \ref{pf_constructionoftheconjugacy}) satisfies closeness estimates optimal in $h$ near the PF point: there exists
a positive constant $const>0$ such that for all $0<h\le h_0$, $|x|\le \varepsilon_0$ and $|\alpha|\le \alpha_0$ \[|x-J(h,x,\alpha)|\le const\cdot h^p.\]
\end{thm}

Similarly to the TC case, in some regions we have $\textit{const}\cdot
 \alpha^\kappa\,h^p$ on the right-hand side with a suitable $\kappa=1$ or $\kappa=2$.

\section{Conjugacy in the Discretized TC Bifurcation}\label{conjugacyTC}

\subsection{Construction of the conjugacy in the TC case}\label{tk_constructionoftheconjugacy}

According to Section \ref{TCmainresultssection}, it is enough to define a conjugacy
map between the corresponding normal forms, and these normal forms 
can be chosen as
\begin{equation}\label{tk_normalFi}
\NFi(h,x,\al):=(1+h\al)x+ h x^2+h x^3\,\eta_3(h,x,\al)
\end{equation}
and
\begin{equation}\label{tk_normalfi}
\Nfi(h,x,\al):=(1+h\al)x+ h x^2+h x^3\,{\widetilde{\eta}_3}(h,x,\al),
\end{equation}
where $\eta_3$ and $\widetilde{\eta}_3$ are smooth functions.  
Let $K>0$ denote a uniform bound on $\Big|\frac{\dd^i}{\dx^i}\,\eta(h,\cdot,\al)\Big|$ ($i\in\{0,1,2\}$,\,$\eta\in\{\eta_3,\widetilde{\eta}_3\}$) in a neighborhood of the origin for
small $h>0$ and $|\al|$, as well as a uniform bound on $\Big|\frac{\dd}{\dd\al}\,\eta(h,x,\cdot)\Big|$ ($\eta\in\{\eta_3,\widetilde{\eta}_3\}$) in a neighborhood of the origin for 
small $h>0$ and $|x|$. 
Theorem \ref{phinormalforma and closeness_tk} and the reductions mentioned afterwards imply
that there exists a constant $c>0$ such that 
\begin{equation}\label{tk_normalformakkozelsege}
\left| \NFi(h,x,\al)-\Nfi(h,x,\al)\right|\le c\cdot h^{p+1}|x|^3
\end{equation}
holds for all sufficiently small $h>0$, $|x|\ge 0$ and $|\al|\ge 0$. Throughout Section
\ref{conjugacyTC}, $c$ will denote this particular
positive constant. Other generic positive constants, if needed, are denoted by $const$.

\begin{figure}[h]
 \centering
 \includegraphics{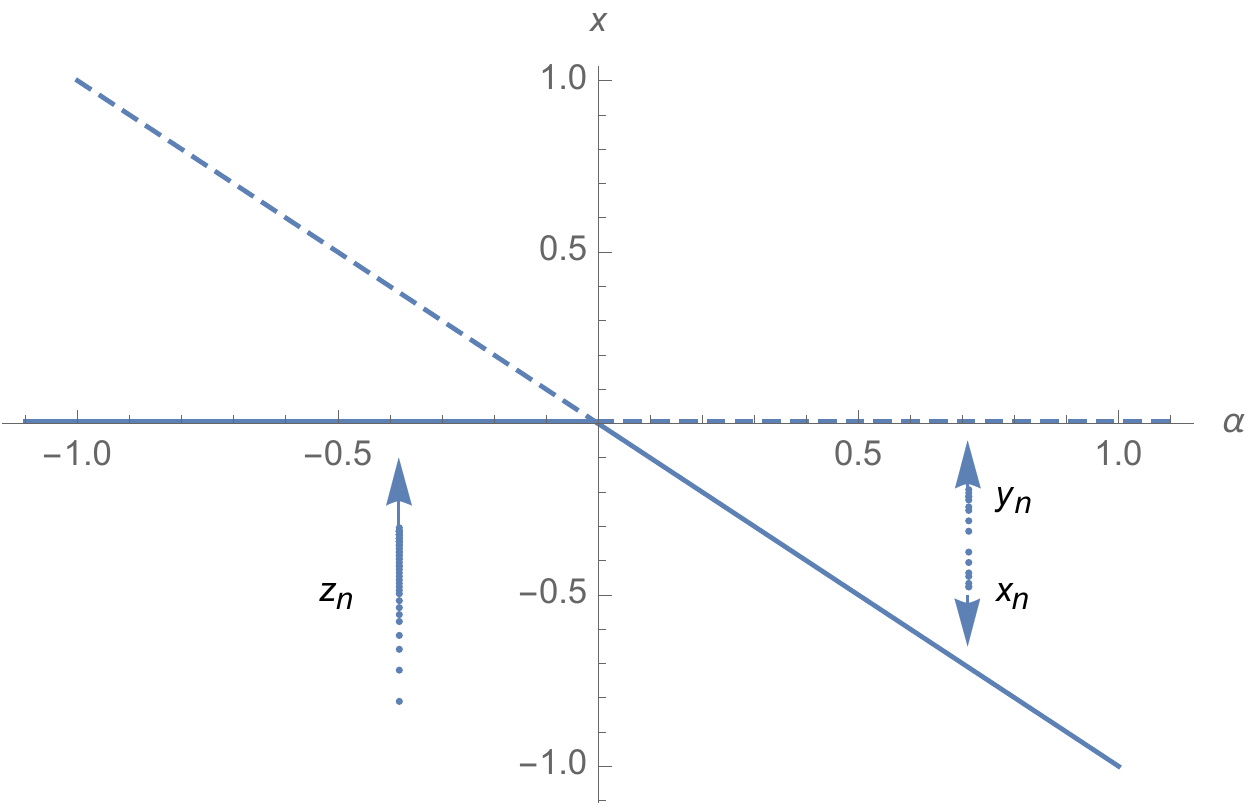}
 \caption{The branch of stable (continuous line) and unstable (dashed line) fixed points of $\Nfie$ 
are shown together with the first few terms of the inner sequences
$x_n(h,\alpha)$ and $y_n(h,\alpha)$, and the outer sequence $z_n(h,\alpha)$, for some fixed $h>0$ and $\al$. The arrows indicate the direction of these sequences.}\label{figure3}
\end{figure}

It is easy to see that for any fixed $h>0$,
$\omFinul(h,\al) = 0$ is an attracting fixed point of the map $\NFi(h,\cdot,\al)$ for $\al<0$,
and repelling for $\al>0$, see Figure \ref{figure3}. For any fixed $\al\in [-\al_0,\al_0]\setminus \{0\}$, $\NFi(h,\cdot,\al)$ 
possesses another fixed point, denoted by $\omFiplusz\equiv\omFiplusz(h,\al)>0$ (for $\al<0$) and $\omFi\equiv\omFi(h,\al)<0$ (for $\al>0$). It is seen that $\omFiplusz$ is repelling and $\omFi$ is attracting. The two branches of fixed points, $\omFinul(h,\al)$ and $\omega_{\Phi,_{\pm}}(h,\al)$ merge at $\al=0$. Analogous results hold for the map $\Nfi(h,\cdot,\al)$. Its fixed points are denoted by $\omfinul \equiv 0$, $\omfi$ and $\omfiplusz$.

In the rest of this section
the conjugacy map is constructed in a natural way in the $x\le 0$ region.
The $x>0$ case is similar due to symmetry, hence it is omitted here. In what follows, we suppose that 
\[
0<h\le h_0:=\frac{1}{5},
\]
\begin{equation}\label{tk_0kicsi}
 |x|\le \varepsilon_0:=\min\left(\frac{1}{25},\frac{1}{25K}\right)\ \mathrm{and\ } 
\end{equation} 
\[
|\al| \le \al_0:=\min\left(\frac{1}{51},\frac{1}{51K}\right).
\]
With these values of $h_0$, $\varepsilon_0$ and $\al_0$, all constructions and proofs below can be carried out. (There is only one constraint which has not been taken into account explicitly: if the domain of definition of the functions $\eta_3$ and $\widetilde{\eta}_3$ is smaller than $(0,h_0]\times [-\varepsilon_0,\varepsilon_0]\times [-\al_0,\al_0]$ given above, then $h_0$, $\varepsilon_0$ or $\al_0$ should be decreased further suitably.)

By iterating one of the normal forms, say $\Nfi(h,\cdot,\al)$, let us define three sequences $x_n$, $y_n$  and $z_n$. For $\al>0$, let $x_n\equiv x_n(h,\al)$ be defined as
\[
x_{n+1}:={\cal{N}_{\vfi}}(h,x_n,\al), \quad n=0,1,2,\ldots
\]
with $x_0:=-\frac{\al}{3}$, and let $y_n\equiv y_n(h,\al)$ be defined as
\[
y_{n}:=\left(\Nfi^E\right)^{[-n]}(x_0), \quad n=0,1,2,\ldots,
\]
so $y_0:=x_0$, and set $y_{-1}:=x_1$. Finally, for all $\al\in [-\al_0,\al_0]$ define $z_n\equiv z_n(h,\al)$ as
\[
z_{n}:=\left(\Nfi^E\right)^{[n]}(z_0), \quad n=0,1,2,\ldots,
\]
with $z_0<0$ being independent of $h$ and $\al$ such that $2\al_0<|z_0|<\frac{1}{2K}$ holds. An appropriate choice for $z_0$ is, e.g., $z_0:=-\varepsilon_0$.

Simple calculations show that under conditions (\ref{tk_0kicsi}), for example, both $\Nfi^E$ and $\NFi^E$ (together with their inverses) are strictly monotone increasing, moreover, $|\al|<\frac{6}{K}$ implies  $x_0(\al)>x_1(h,\al)$, and $2\al_0<|z_0|<\frac{1}{2K}$ implies $z_0<z_1(h,\al)$.  
This means that $x_n$ is strictly decreasing, $y_n$ is strictly increasing (for $\al>0$ and $n\ge 0$), and $\lim_{n\rightarrow \infty}x_n(h,\al)=\omfi$, while $\lim_{n\rightarrow \infty}y_n(h,\al)=\omfinul$.
Similarly, $z_n$ is strictly increasing, moreover,  $\lim_{n\rightarrow \infty}z_n(h,\al)=\omfi$
for $\al>0$, and $\lim_{n\rightarrow \infty}z_n(h,\al)=\omfinul$ for $\al\le 0$.

 A homeomorphism $J^E$ satisfying the conjugacy equation \begin{equation}\label{tk_konjugacios alapegyenlet} J^E\circ \Nfie=\NFi^E\circ J^E\end{equation} is now piecewise defined on the fundamental domains, i.e., on  $[x_{n+1},x_{n}]$, $[y_{n},y_{n+1}]$ and $[z_{n},z_{n+1}]$ ($n\in\mathbb{N}$), 
for any fixed $0<h\le h_0$ and $-\al_0\le \al \le \al_0$.

We first consider the inner region, that is, the region between the fixed points for $0<\al\le\al_0$.
Set $J^E(x_0):=x_0$, $J^E(x_1):=\NFi^E(x_0)$ and for $x\in [x_1,x_0]$ extend $J^E$ linearly. For $n\ge 1$ and
$x\in [x_{n+1},x_n]$, we recursively 
set 
\[
J^E(x):=\left(\NFi^E\circ J^E \circ \left(\Nfi^E\right)^{[-1]}\right)(x),
\]
while for $n\ge 0$ and $x\in [y_n,y_{n+1}]$, we let
\[
J^E(x):=\left(\left(\NFi^E\right)^{[-1]} \circ J^E \circ \Nfie\right)(x).
\]
(Since $[y_{-1},y_0] \equiv [x_1,x_0]$, these two definitions are compatible.) Finally, set
$J^E(\omfi):=\omFi$ and $J^E(\omfinul):=\omFinul$. Then $J^E$ is continuous, strictly monotone increasing on
$[\omfi,0]$, since it is a composition of three such functions, and satisfies
(\ref{tk_konjugacios alapegyenlet}).

In the outer region, i.e., below the fixed points, fix $z_0<0$ ($2\al_0<|z_0|<\frac{1}{2K}$).
For $\al\in [-\al_0,\al_0]$, the definition of $J^E$ is
analogous to the above construction  with  $x_n$: this time $z_n$ plays the role of $x_n$ (now the counterpart of the sequence $y_n$ is not needed).
Then the function $J^E$ is continuous and strictly monotone increasing on $[z_0,\omfi]$ 
(for $0<\al\le \al_0$) and on $[z_0,\omfinul]$ (for $-\al_0\le \al\le 0$), and satisfies (\ref{tk_konjugacios alapegyenlet}).

The construction of $J^E$---with the appropriate and natural modifications---in the upper half-plane $x>0$ is analogous to the one presented above.

\subsection{Conjugacy estimates: optimality at the TC fixed points}\label{tk_optimality}

We  prove here that the constructed conjugacy $J^E$ is $\cal{O}$$\left(h^p\,\al^2\right)$-close to the
identity at the fixed points $\omfi(h,\al)$, moreover, an explicit example will show
that this estimate is optimal in $h$ and $\al$. Since fixed points must be mapped
into nearby fixed points by the conjugacy,
this means that the estimates in Theorem
\ref{Conjugacy estimate in the TC bifurcation} are optimal in $h$.

First, two useful lemmas are presented.
\begin{lemma}\label{tk_Nfi fixpontja}
For every $0<h\le h_0$ and $0< \al \le \al_0$ we have that
\[
\{\omfi,\omFi\} \subset \left(-\frac{3}{2}\al,-\frac{6}{7}\al\right ).
\]
\end{lemma}
\textit{Proof.} By definition, $\omfi$ solves
$\al+x+x^2\cdot\widetilde{\eta}_3(h,x,\al)=0$. But $|x|\le\frac{1}{6K}$ implies $\frac{2}{3}\le 1+x\,\widetilde{\eta}_3\le \frac{7}{6}$, so \[
-\frac{3\al}{2}\le\omfi=\frac{-\al}{1+\omfi\cdot\widetilde{\eta}_3(h,\omfi,\al)} \le
-\frac{6\al}{7}.
\]
The proof for $\omFi$ is similar.\quad $\Box$

\begin{lemma}\label{tk_NFidervialtbecsles} For any $0<h\le h_0$, $-\varepsilon_0\le x<0$ and $-\al_0\le \al\le \al_0$, we have that
\[
\left(\NFi^E\right)^{\prime}(x)\le 1+h\al+\frac{7}{4}h x.
\]
\end{lemma}
\textit{Proof.} The conditions in (\ref{tk_0kicsi}) have been set up to imply this inequality, too.
\quad $\Box$

\begin{lemma}\label{tk_fixpontoktavolsaga}
 For any $0<h\le h_0$ and $0<\al\le \al_0$ satisfying (\ref{tk_0kicsi}), we have that
\[
|\,\omfi-\omFi|\le \frac{27}{4}c\cdot h^{p}\,\al^2.
\]
\end{lemma}
\textit{Proof.} By Lemma \ref{tk_Nfi fixpontja}, (\ref{tk_normalformakkozelsege}) and Lemma \ref{tk_NFidervialtbecsles}, we have
\[|id-J^E|(\omfi(h,\al))\le |\Nfie(\omfi)-\NFi^E(\omfi)|+ |\NFi^E(\omfi)-\NFi^E(\omFi)|\le 
\]
\[
c\cdot h^{p+1}|\omfi|^3+\left(\sup_{[\{\,\omfi,\omFi\}]}\left(\NFi^E\right)^{\prime}\right)|\,\omfi-\omFi|\le
\]
\[
\frac{27}{8}c\cdot h^{p+1}\al^3+\left(1-\frac{h\al}{2}\right)|\,\omfi-\omFi|.
\]
Solving the above inequality for $|\,\omfi-\omFi|\equiv |id-J^E|(\omfi)$ yields the desired result.
\quad $\Box$

\begin{rem} The following
example shows that the distance of fixed points of normal forms 
satisfying (\ref{tk_normalformakkozelsege}) can be bounded from below by
$const\cdot h^p$ ($h\rightarrow 0$). 
Indeed, set 
$\NFi(h,x,\al):=(1+h\al)x+h x^2$ and $\Nfi(h,x,\al):=(1+h\al)x+h x^2+h^{p+1}x^3$. Then these maps satisfy
(\ref{tk_normalformakkozelsege}) in a neighborhood of the origin, moreover, $\omFi=-\al$ and 
$\omfi=\frac{-1+\sqrt{1-4h^p\,\al}}{2h^p}$.
Using inequality $1+\frac{t}{2}-\frac{t^2}{4}\le\sqrt{1+t}\le 1+\frac{t}{2}-\frac{t^2}{8}$ for $-\frac{1}{2}\le t\le 0$, one sees that
\[
|\,\omfi-\omFi|\ge h^p\, \al^2,
\]
if, for example, $h\le 1$ and $\al\le \frac{1}{8}$.\quad $\Box$
\end{rem}

\subsection{Conjugacy estimates in the inner TC region}\label{tk_innerregion}

In this section, for any fixed $0<h\le h_0$ and $0<\al\le \al_0$, the closeness estimate is proved in $(\omfi,0)$. First, the proof is presented in $(\omfi,x_0]$ (Lemma \ref{tk_inner_1}),
then in $(x_0,0)$ (Lemma \ref{tk_lemma3.11inner2}). 
(As for the boundary points, the conjugacy estimates at $\omfi$ have already been covered by Lemma 
\ref{tk_fixpontoktavolsaga}, while at $0$ these estimates are trivial.)

It is clear that $\sup_{(\omfi,x_0]}|id-J^E| =  \sup_{n\in \mathbb{N}} \sup_{[x_{n+1},x_{n}]}|id-J^E|$.
Since $x_0=J^E(x_0)$, we have that
\[
\sup_{[x_1,x_{0}]}|id-J^E|= |x_1-J^E(x_1)|=|\Nfi^E(x_0)-\NFi^E(x_0)|\le
\]
\[
c\cdot h^{p+1} |x_0|^3=\frac{c}{27}h^{p+1}\al^3,
\]
while for $n\ge 1$,
\[
\sup_{[x_{n+1},x_n]}|id-J^E|\le  \sup_{[x_{n+1},x_n]}\left|\Nfi^E\circ (\Nfi^E)^{[-1]}-\NFi^E\circ  (\Nfi^E)^{[-1]}
\right|+ \] \[ \sup_{[x_{n+1},x_n]}\left|\NFi^E\circ (\Nfi^E)^{[-1]}-\NFi^E\circ J^E\circ (\Nfi^E)^{[-1]}\right|= \]
\[\sup_{[x_n,x_{n-1}]}\left|\Nfi^E-\NFi^E\right|+
\sup_{[x_n,x_{n-1}]}\left|\NFi^E-\NFi^E\circ J^E\right| \le \]
\[\sup_{[x_n,x_{n-1}]}\left|\Nfi^E-\NFi^E\right|+ \sup_{x\in [x_n,x_{n-1}]} \left(\left(\sup_{[\{x,J^E(x)\}]}
(\NFi^E)^{\prime}\right)|\,x-J^E(x)|\right) \le \]
\[
c\cdot h^{p+1}  |x_n|^3+ \left(1+h\al+\frac{7}{4}h\max\left(x_{n-1},J^E(x_{n-1})\right)\right) \sup_{[x_n,x_{n-1}]}|id-J^E|.
\]
Here the last inequality is true because we take into account Lemma \ref{tk_NFidervialtbecsles} and use the fact that the functions $id$ and $J^E$ are increasing, so
\[
\sup_{[\{x,J^E(x)\}]}(\NFi^E)^{\prime}\le \sup_{[\{x,J^E(x)\}]}(1+h\al+\frac{7}{4}h\cdot id)\le 1+h\al+\frac{7}{4}h\max\left(x,J^E(x)\right).
\]

From these we have for $n\ge 1$ that
\[
\sup_{[x_{n+1},x_n]}|id-J^E|\le c\cdot h^{p+1} \sum_{i=0}^{n}|x_i|^3\prod_{j=i}^{n-1}\left(1+h\al+\frac{7}{4}h\max\left(x_{j},J^E(x_{j})\right)\right),
\]
where $\prod_{j=n}^{n-1}$ is understood to be $1$.
 
So in order to prove that the conjugacy $J^E$ is $\cal{O}$$(h^p)$-close to the identity
on the interval $(\omfi,x_0]$ for any $h\in (0,h_0]$ and $\al \in (0,\al_0]$, it is enough to show that
\begin{equation}\label{tk_innerbecslese}
\sup_{h\in (0,h_0]}\sup_{\al\in (0,\al_0]}\sup_{n\in \mathbb{N}} h \sum_{i=0}^{n} 
|x_i|^3\prod_{j=i}^{n-1}\left(1+h\al+\frac{7}{4}h\max\left(x_{j},J^E(x_{j})\right)\right)\le const
\end{equation}
holds with a suitable $const\ge 0$.

First an explicit estimate of the sequence $\max\left(x_{n},J^E(x_{n})\right)$ is given. 

\begin{lemma}\label{tk_xknagysagrend}
For $n\ge 0$, set \[a_n(h,\al):=-\frac{3}{4}\al\cdot\frac{(1+h\al)^{n+1}}{2+(1+h\al)^n}.\] 
Then $x_n\in (\omfi,a_n)$ and $J^E(x_n)\in (\omFi,a_n)$.
\end{lemma}
\textit{Proof.} It is easily checked that, due to assumptions (\ref{tk_0kicsi}), 
\[\max\left(\omfi,\omFi\right)<a_n\] for $n\ge 0$, so the intervals
in the lemma are non-degenerate. We proceed by induction:
$a_0=-\frac{\al}{4}(1+h\al)>x_0\equiv J^E(x_0)\equiv -\frac{\al}{3}$ is equivalent to $h\al<\frac{1}{3}$, being true by assumptions (\ref{tk_0kicsi}) on $h_0$ and $\al_0$. 
So suppose that the statement is true for some $n\ge 0$. Since $\Nfie(x)<(1+h\al)x+\frac{6}{5}h x^2$ is implied by $|x|\le \varepsilon_0 <\frac{1}{5K}$, and $\Nfi^E$ is monotone increasing, we get that
\[
x_{n+1}=\Nfie(x_n)<\Nfie(a_n)<(1+h\al)a_n+\frac{6}{5}h a_n^2.
\]
Therefore it is enough to prove that the right-hand side above is smaller than $a_{n+1}$. But 
\[
a_{n+1}-\left((1+h\al)a_n+\frac{6}{5}h a_n^2\right)=
\]
\[
-\frac{3h\al^2(1+h\al)^{2+2n}\left(-2+(1+h\al)^n(-1+9h\al)\right)}{40\left(2+(1+h\al)^n\right)^2\left(2+(1+h\al)^{n+1}\right)}>0
\]
is equivalent to $-2+(1+h\al)^n(-1+9h\al)<0$, which is implied by $h\al<\frac{1}{9}$.

Of course, the above inequalities remain true, if $\Nfi$ is replaced by $\NFi$, also noticing that, by construction, $J^E(x_{n+1})=\NFi^E(J^E(x_n))$, so the proof is complete.
\quad $\Box$

\begin{rem}
The induction would fail, if, in estimate $\Nfie(x)<(1+h\al)x+\frac{6}{5}h x^2$, the constant $\frac{6}{5}$ was replaced by, say, $\frac{7}{5}$. (The explanation resides in the particular choice of the constant $\frac{3}{4}$ in the definition of $a_n$, since $\frac{3}{4}\cdot\frac{6}{5}<1<
\frac{3}{4}\cdot\frac{7}{5}$.)\quad $\Box$
\end{rem}

\begin{rem} The upper estimate $a_n$ has been found by experiments with \textit{Mathematica} based on the parameterized model functions 
(\ref{thulsmodelfunction}).\quad $\Box$
\end{rem} 

In order to prove the boundedness of the left-hand side of (\ref{tk_innerbecslese}), the sum $\sum_{i=0}^n$ will
be split into two. An appropriate index to split at is $\lceil \frac{const}{h\al}\rceil$, as established by the following lemma.

\begin{lemma}\label{tk_inner split}
Suppose that $n>\lceil \frac{6}{h\al}\rceil$. Then
\[
\max\left(x_n,J^E(x_n)\right)<-\frac{2}{3}\al,
\]
hence 
\[
1+h\al+\frac{7}{4}h\max\left(x_n,J^E(x_n)\right)<1-\frac{h\al}{6}
\]
holds for $n>\lceil \frac{6}{h\al}\rceil$.
\end{lemma}
\textit{Proof.} By Lemma \ref{tk_xknagysagrend} it is sufficient to show that
$n>\lceil \frac{6}{h\al}\rceil$ implies $a_n<-\frac{2}{3}\al$. This latter inequality is equivalent to $(1+h\al)^n (1+9h\al)>16$. But for $n>\lceil \frac{6}{h\al}\rceil$, 
\[
(1+h\al)^n>\left(1+h\al\right)^{\lceil \frac{6}{h\al}\rceil}=\left(1+\frac{1}{\frac{1}{h\al}}\right)^{\left(1+\frac{1}{h\al}\right)\cdot\frac{h\al}{1+h\al}\cdot \lceil \frac{6}{h\al}\rceil}.
\]
However, it is known that $\left(1+\frac{1}{A}\right)^{A+1}>e$ for $A\ge 1$, and it is 
easy to see that 
$\frac{B}{1+B}\cdot \lceil \frac{6}{B}\rceil>3$ for $0<B<1$. Since 
$e^3>16$, the proof is complete.
\quad $\Box$

Now we can turn to (\ref{tk_innerbecslese}). Fix $h\in (0,h_0]$, $\al\in (0,\al_0]$ and $n\in \mathbb{N}^+$. (For 
$n\le \lceil \frac{6}{h\al}\rceil$,  the sums $\sum_{i=\lceil \frac{6}{h\al}\rceil+1}^n$ below are, of course, not present, making the proof even simpler.) Since now $\omfi<x_i<0$, by Lemma \ref{tk_Nfi fixpontja} we have $|x_i|\le \frac{3}{2}\al$, and by monotonicity $\max\left(x_j,J^E(x_j)\right)\le x_0\equiv J^E(x_0)\equiv -\frac{\al}{3}$. Moreover, by using Lemma \ref{tk_inner split}, assumption $h\al<1$ from (\ref{tk_0kicsi}), and inequality $(1+\frac{1}{A})^A\le e$ (for $A\ge 1$), we get that
\[
h \sum_{i=0}^{n} 
|x_i|^3\prod_{j=i}^{n-1}\left(1+h\al+\frac{7}{4}h\max\left(x_{j},J^E(x_{j})\right)\right)\le 
\]
\[
\frac{27h\al^3}{8}\sum_{i=0}^{\lceil \frac{6}{h\al}\rceil} 
\prod_{j=1}^{\lceil \frac{6}{h\al}\rceil-1}\left(1+h\al-\frac{7}{4}\cdot\frac{h\al}{3}\right) 
+\frac{27h\al^3}{8}\sum_{i=\lceil \frac{6}{h\al}\rceil+1}^{n}\prod_{j=i}^{n-1}\left(1-\frac{h\al}{6}\right)\le 
\]
\[
\frac{27h\al^3}{8}\left(1+\frac{5}{12}h\al\right)^{\frac{6}{h\al}}\left(\Big{\lceil}{\frac{6}{h\al}}\Big{\rceil}+1\right)+
\frac{27h\al^3}{8}\sum_{i=\lceil \frac{6}{h\al}\rceil+1}^{n}\left(1-\frac{h\al}{6}\right)^{n-i}\le
\]
\[
\frac{27h\al^3}{8}\left(1+\frac{5}{12}h\al\right)^{\frac{12}{5h\al}\cdot\frac{5h\al}{12}\cdot\frac{6}{h\al}}\left(\frac{6+2h\al}{h\al}\right)+
\frac{27h\al^3}{8}\sum_{i=0}^{\infty}\left(1-\frac{h\al}{6}\right)^{i}\le
\]
\[
\frac{27h\al^3}{8}\cdot e^{\frac{30}{12}}\cdot \frac{8}{h\al}+
\frac{27h\al^3}{8}\cdot\frac{6}{h\al}\le 350\,\al^2.
\]

\noindent Therefore, $\sup_{[x_{n+1},x_{n}]}|id-J^E|\le 350c\cdot h^p \al^2$  for any $h\in (0,h_0]$, $\al\in (0,\al_0]$ and $n\ge 1$, moreover, as we have seen, $\sup_{[x_{1},x_{0}]}|id-J^E|\le \frac{c}{27}h^{p+1}\al^3$, which yield the following lemma.
\begin{lemma}\label{tk_inner_1}
Under assumption (\ref{tk_0kicsi}),
\[
\sup_{(\omfi,x_0]}|id-J^E|\le 350c\cdot h^p \al^2.
\]
\end{lemma}
In the rest of this section, the closeness estimate is proved in the interval $(y_0,\omfinul)$. Recall that $y_0=x_0=J^E(x_0)\equiv -\frac{\al}{3}$ and $\omfinul=\omFinul\equiv 0$.

Suppose that $n\ge 1$ (the case $n=0$ will be examined later). Then
\[
\sup_{[y_n,y_{n+1}]}|id-J^E|=\sup_{[y_n,y_{n+1}]}
\left|\left(\NFi^E\right)^{[-1]}\circ \NFi^E-\left(\NFi^E\right)^{[-1]}\circ J^E \circ \Nfi^E\right|\le
\] 
\[
\sup_{x\in [y_n,y_{n+1}]}\left[\left(\sup_{[\{\NFi^E(x),J^E\circ \Nfi^E(x)\}]}\left((\NFi^E)^{[-1]}\right)^{\prime}\right) \left( \left|\NFi^E-\Nfi^E\right|(x)+
\left|\Nfi^E-J^E\circ \Nfi^E\right|(x)\right)\right]\le
\]
\[
\left[ \sup_{x\in [y_n,y_{n+1}]}\sup_{[\{\NFi^E(x),J^E\circ \Nfi^E(x)\}]}\left((\NFi^E)^{[-1]}\right)^{\prime}\right]
\left[c\cdot h^{p+1} |y_n|^3+ \sup_{[y_{n-1},y_{n}]} |id-J^E|\right], 
\]
provided that $\sup_{[\{\NFi^E(x),J^E\circ \Nfi^E(x)\}]}\left((\NFi^E)^{[-1]}\right)^{\prime}$ is non-negative.

\begin{lemma}
Suppose that $n\ge 1$, then under assumption (\ref{tk_0kicsi}) we have that
\[
 \sup_{x\in [y_n,y_{n+1}]}\sup_{[\{\NFi^E(x),J^E\circ \Nfi^E(x)\}]}\left((\NFi^E)^{[-1]}\right)^{\prime}\le 1-\frac{h\al}{8}.
\]
\end{lemma}
\textit{Proof.}
\[
\sup_{x\in [y_n,y_{n+1}]}\sup_{[\{\NFi^E(x),J^E\circ \Nfi^E(x)\}]}\left((\NFi^E)^{[-1]}\right)^{\prime}=
\sup_{x\in [y_n,y_{n+1}]}\sup_{[\{\NFi^E(x),J^E\circ \Nfi^E(x)\}]}\frac{1}{(\NFi^E)^{\prime}\circ(\NFi^E)^{[-1]}}=
\]
\[
\sup_{x\in [y_n,y_{n+1}]}\sup_{[\{x,(\NFi^E)^{[-1]}\circ J^E\circ \Nfi^E(x)\}]}\frac{1}{(\NFi^E)^{\prime}}=\ldots
\]
But by definition $(\NFi^E)^{[-1]}\circ J^E\circ \Nfi^E(x)=J^E(x)$ for $x\in [y_n,y_{n+1}]$, and $[\{x,J^E(x)\}]=
[\,\min(x,J^E(x)),\max(x,J^E(x))]$, moreover, by the monotonicity of $id$ and $J^E$ we obtain that 
\[
\ldots = \sup_{[\min(y_n,J^E(y_n)),\max(y_{n+1},J^E(y_{n+1}))]}\frac{1}{(\NFi^E)^{\prime}}\le \ldots
\]
By construction, however, $[\,\min(y_n,J^E(y_n)),\max(y_{n+1},J^E(y_{n+1}))]\subset (y_0,0)=
\left(-\frac{\al}{3},0\right)$ and $\left(\NFi^E\right)^{\prime}$ is nonnegative here by assumption (\ref{tk_0kicsi}), justifying the computations just above the lemma. We now continue the proof of the lemma as
\[
\ldots \le \sup_{\left(-\frac{\al}{3},0\right)}\frac{1}{(\NFi^E)^{\prime}}\le \ldots
\]
It is easy to see that  (\ref{tk_0kicsi}) together with $x<0$ imply  $(\NFi^E)^{\prime}(x)\ge
1+h\al+\frac{9}{4}h x\ge 0$. So
\[
\ldots \le \sup_{x\in \left(-\frac{\al}{3},0\right)}\frac{1}{1+h\al+\frac{9}{4}h x}\le \frac{1}{1+h\al+\frac{9}{4}h \left(-\frac{\al}{3}\right)}=\frac{1}{1+\frac{1}{4}h\al}\le 1-\frac{h\al}{8},
\]
since $\frac{1}{1+A}\le 1-\frac{A}{2}$ for $A\in [0,1]$.
\quad $\Box$

We have thus proved (also using $|y_n|\le \frac{\al}{3}$) that for $n\ge 1$
\begin{equation}\label{tk_[y_n,y_{n+1}]}
\sup_{[y_n,y_{n+1}]}|id-J^E|\le \left(1-\frac{h\al}{8}\right)\left[\frac{c}{27}\cdot h^{p+1} \al^3+ \sup_{[y_{n-1},y_{n}]} |id-J^E|\right].
\end{equation}

For $n=0$, similarly as before, we get that
\[
\sup_{[y_0,y_{1}]}|id-J^E|
\le \left[ \sup_{x\in [y_0,y_{1}]}\sup_{[\{\NFi^E(x),J^E\circ \Nfi^E(x)\}]}\left((\NFi^E)^{[-1]}\right)^{\prime}\right]
\left[c\cdot h^{p+1} |y_0|^3+ \sup_{[y_{-1},y_{0}]} |id-J^E|\right].
\]
But $[y_{-1},y_{0}]\equiv[x_{1},x_{0}]$, so the second factor $[\ldots]$ is bounded by $2\cdot \frac{c}{27}h^{p+1}\al^3$. As for the first factor $[\ldots]$, we notice that $y_0<(\NFi^E)^{[-1]}(y_0)$ (since this is equivalent to $x_1<x_0$), which implies that 
\[
\sup_{x\in [y_0,y_{1}]}\sup_{[\{\NFi^E(x),J^E\circ \Nfi^E(x)\}]}\left((\NFi^E)^{[-1]}\right)^{\prime}
=\sup_{x\in [y_0,y_{1}]}\sup_{[\{x,(\NFi^E)^{[-1]}\circ J^E\circ \Nfi^E(x)\}]}\frac{1}{(\NFi^E)^{\prime}}=
\]  
\[
\sup_{[y_0,y_{1}]\cup [y_0,(\NFi^E)^{[-1]}(y_0)]}\frac{1}{(\NFi^E)^{\prime}}\le
\sup_{[y_0,0)}\frac{1}{(\NFi^E)^{\prime}}\le 1.
\]
Therefore
\begin{equation}\label{tk_[y_0,y_{1}]}
\sup_{[y_0,y_{1}]}|id-J^E|\le 2\cdot \frac{c}{27}h^{p+1}\al^3.
\end{equation}
\indent Repeated application of (\ref{tk_[y_n,y_{n+1}]}), moreover (\ref{tk_[y_0,y_{1}]}) yield for $n\ge 1$ that
\[
\sup_{[y_n,y_{n+1}]}|id-J^E|\le \left(1-\frac{h\al}{8} \right)^n \sup_{[y_0,y_{1}]}|id-J^E| + \frac{c}{27}h^{p+1}\al^3\sum_{i=1}^{n} \left(1-\frac{h\al}{8} \right)^i\le
\]
\[
1\cdot 2\cdot \frac{c}{27}h^{p+1}\al^3+\frac{c}{27}h^{p+1}\al^3 \cdot \frac{8}{h\al}\le \frac{c}{3}h^{p}\al^2, 
\]
due to $h\al\le \frac{1}{2}$ by  (\ref{tk_0kicsi}). The same upper estimate is valid for $n=0$, so we have proved the following result.
\begin{lemma}\label{tk_lemma3.11inner2}
Under assumption (\ref{tk_0kicsi})
\[
\sup_{(x_0,0)}|id-J^E|\le \frac{c}{3}h^{p}\al^2.
\]
\end{lemma}

\subsection{Conjugacy estimates in the outer TC region}\label{tk_outerregion}

In this section, we first prove an ${\cal{O}}(h^p)$-closeness estimate in the interval $[z_0,\omfi)$ for $\al>0$, then, in the second part, the closeness is proved on $[z_0,\omFinul)\equiv [z_0,0)$ for $\al\le 0$.

The derivation of the following formulae is similar to their counterparts in the
inner region, with the difference that now an extra term and an index-shift occur
(since this time the sequence $z_n$ is increasing).
For $n\ge 1$ (by also using $(\ref{tk_0kicsi})$) we have that
\[
\sup_{[z_n,z_{n+1}]}|id-J^E|\le c\cdot h^{p+1} |z_0|^3\prod_{j=1}^{n}\left(1+h\al+\frac{7}{4}h\max\left(z_{j},J^E(z_{j})\right)\right)+
\]
\begin{equation}\label{tk_outer}
c\cdot h^{p+1} \sum_{i=0}^{n-1}|z_i|^3\prod_{j=i+2}^{n}\left(1+h\al+\frac{7}{4}h\max\left(z_{j},J^E(z_{j})\right)\right),
\end{equation}
where, again, $\prod_{j={n+1}}^{n}$ is 1, while the $n=0$ case is simply 
\[
\sup_{[z_0,z_{1}]}|id-J^E|\le c\cdot h^{p+1} |z_0|^3.
\]

The following main lemma, as a counterpart of Lemma \ref{tk_xknagysagrend}, gives a lower estimate of the
sequence $z_n$ for $\al>0$.

\begin{lemma}\label{tk_zknagysagrend}
For $n\ge 0$, set \[b_n(h,\al):=-2\al\cdot\frac{(1+h\al)^{n+1}}{-1+\al+(1+h\al)^n},\] then 
$b_n\le \min\left(z_n,J^E(z_n)\right)$.
\end{lemma}
\textit{Proof.}
We see that $b_0=-2-2h\al<-2\le -1\le -\varepsilon_0\le z_0=J^E(z_0)$ holds due to assumption (\ref{tk_0kicsi}). So suppose that the statement is true for some $n\ge 0$. Since $\Nfie(x)\ge (1+h\al)x+\frac{3}{5}h x^2$ follows from $|x|\le \varepsilon_0<\frac{2}{5K}$, moreover, $(1+h\al)id+\frac{3}{5}h\,id^2$ is monotone increasing (which is implied by, e.g., $|x|\le \frac{5}{6h}$; but it is easy to see that $h\le\frac{5}{18}$ and $-3<b_n<0$ follows from  (\ref{tk_0kicsi}), hence $|b_n|\le \frac{5}{6h}$), so we obtain that
\[
z_{n+1}=\Nfie(z_n)\ge (1+h\al)z_n+\frac{3}{5}h z_n^2\ge (1+h\al)b_n+\frac{3}{5}h b_n^2.
\]
Therefore it is sufficient to show that 
$(1+h\al)b_n+\frac{3}{5}h b_n^2\ge b_{n+1}$.
However, the last inequality is equivalent to
\[
0\le \frac{2h\al^2 (1+h\al)^{2+2n}}{5\left(-1+\al+(1+h\al)^n\right)^2}\cdot\frac{-1+\al+(1+h\al)^n(1+6h\al)}{-1+\al+(1+h\al)^{n+1}}\,,
\]
which is true since $\al>0$ and $h>0$.

The proof remains valid if $\Nfi$ is replaced by $\NFi$ (and $J^E(z_n)$ is written instead of $z_n$), hence $b_n\le J^E(z_n)$ also holds.
\quad $\Box$

Now, since $z_j<\omfi$ and $J^E(z_j)<\omFi$, we get  by Lemma \ref{tk_Nfi fixpontja} that the 
right-hand side of (\ref{tk_outer}) is at most 
\[
c\cdot h^{p+1} |z_0|^3\prod_{j=1}^{n}\left(1-\frac{h\al}{2}\right)+ c\cdot h^{p+1} \sum_{i=0}^{n-1}|z_i|^3\prod_{j=i+2}^{n}\left(1-\frac{h\al}{2}\right)\le
\]
\[
c\cdot h^{p+1} |z_0|^3 + c\cdot h^{p+1} \sum_{i=0}^{n-1}|z_i|^3\left(1-\frac{h\al}{2}\right)^{n-1-i}.
\]

We will verify that $h\sum_{i=0}^{n}|z_i|^3\left(1-\frac{h\al}{2}\right)^{n-i}$ is uniformly bounded for any $n\ge 0$,
$0<h\le h_0$ and $0<\al\le \al_0$.

For $n\ge\lceil\frac{1}{h\al}\rceil$,  Lemma \ref{tk_zknagysagrend} says (also using $h\al\le \frac{1}{9}$ and $z_j<0$) that
\[
h\sum_{i=\lceil\frac{1}{h\al}\rceil}^{n}|z_i|^3\left(1-\frac{h\al}{2}\right)^{n-i}\le
h\sum_{i=\lceil\frac{1}{h\al}\rceil}^{n}|b_i|^3\left(1-\frac{h\al}{2}\right)^{n-i}\le 
\]
\[
11 h\al^3\sum_{i=\lceil\frac{1}{h\al}\rceil}^{n}\left(\frac{(1+h\al)^{i}}{-1+\al+(1+h\al)^i}\right)^3\left(1-\frac{h\al}{2}\right)^{n-i}\le \ldots
\]
For these $i$ indices however $\frac{(1+h\al)^{i}}{-1+\al+(1+h\al)^i}\le 3$ holds (since this is implied by 
$\frac{3}{2}\le (1+h\al)^i$, being true by $(1+h\al)^i\ge (1+h\al)^{\frac{1}{h\al}}\ge 1+\frac{1}{h\al}\cdot h\al>\frac{3}{2}$), thus 
\[
\ldots \le 27\cdot 11 \al^2 h\al\sum_{i=0}^{\infty}\left(1-\frac{h\al}{2}\right)^{i} = 594\al^2.
\]
On the other hand, by using that $|z_i|\le 1$ and $h\al\le \frac{1}{9}$ again, 
for $n<\lceil\frac{1}{h\al}\rceil$  we have
\begin{equation}\label{tk_kitevocsokkentes}
h\sum_{i=0}^{n}|z_i|^3\left(1-\frac{h\al}{2}\right)^{n-i}\le 
h\sum_{i=0}^{n}|z_i|^2\left(1-\frac{h\al}{2}\right)^{n-i} \le 
\end{equation}
\[
5 h\sum_{i=0}^{n}\left(\frac{\al(1+h\al)^{i}}{-1+\al+(1+h\al)^i}\right)^2\left(1-\frac{h\al}{2}\right)^{n-i}\le \ldots
\]
Now by the inequalities $e^{\frac{x}{2}}\le 1+x$ (for $x\in [0,1]$) and $1+x\le e^x$ (for $x\in \mathbb{R}$) we get that $(1+h\al)^{2i}\le e^{h\al 2i}\le e^{h\al 2n}\le e^2<8$, 
moreover $\left(1-\frac{h\al}{2}\right)^{n-i}\le e^{-\frac{h\al}{2}(n-i)}$ and $e^{\frac{h\al}{2}i}\le (1+h\al)^{i}$, 
therefore
\begin{equation}\label{40-es_tk}
\ldots \le 40 h\sum_{i=0}^{n}\left(\frac{\al e^{-\frac{h\al}{4}(n-i)}}{-1+\al+e^{\frac{h\al}{2}i}}\right)^2. 
\end{equation}
Set $g_{h,\al,n}(x)\equiv g(x):=\left(\frac{\al \exp\left({-\frac{1}{4}h\al(n-x)}\right)}{-1+\al+\exp\left({\frac{1}{2}h\al x}\right)}\right)^2$ for $x\in [0,\infty)$. Notice that $g$ is bounded at $x=0$. For this function we have that
\[
g^\prime(x)=-\frac{1}{2}h\al^3\, e^{-\frac{1}{2}h\al ( n - x)}\cdot\frac{ 1 - \al + 
        e^{\frac{1}{2}h x \al}} { \left(-1 + \al + e^{\frac{1}{2}h x \al} \right)^3}\, ,
\]
meaning that $g$ is strictly monotone decreasing for $\al<1$. Hence
\[
40 h\sum_{i=0}^{n}\left(\frac{\al e^{-\frac{h\al}{4}(n-i)}}{-1+\al+e^{\frac{h\al}{2}i}}\right)^2=40h+40h\sum_{i=1}^{n}g_{h,\al,n}(i)\le
\]
\[ 
40h + 40 h\int_{0}^{n}g_{h,\al,n}(x)\dx=
40h + 40 h\left[
-2\al\frac{\exp\left({-\frac{1}{2}h\al n}\right)}
{h\left(
-1+\al+\exp\left(\frac{1}{2}h\al x\right)
\right)}
\right]_{x=0}^{n}=
\]
\[
40h+40h\left(
\frac{2\left(1-\exp\left(-\frac{1}{2}h\al n\right)\right)}
{h\left(\exp\left(\frac{1}{2}h\al n\right)-1+\al\right)}
\right)\le 40h+80\left(
\frac{1-\exp\left(-\frac{1}{2}h\al n\right)}
{\exp\left(\frac{1}{2}h\al n\right)-1}
\right)=
\]
\[
40h+80 e^{-\frac{1}{2}h\al n}\le 120,
\]
since $h\le 1$.

Now combining all the estimates so far in the section, under assumption (\ref{tk_0kicsi}) we get
  for $\al>0$ that 
\[
\sup_{[z_0,\omfi)}|id-J^E|=  \sup_{n\in \mathbb{N}} \sup_{[z_{n},z_{n+1}]}|id-J^E|\le
\]
\[
\sup_{n\in \mathbb{N}} \max\left(c\cdot h^{p+1} |z_0|^3,\,c\cdot h^{p+1} |z_0|^3 + c\cdot h^{p+1} \sum_{i=0}^{n}|z_i|^3\left(1-\frac{h\al}{2}\right)^{n-i}\right)\le
\]
\[
c\cdot h^{p+1} |z_0|^3 + c\cdot h^{p}\cdot (120+594\al^2)\le 130c\cdot h^p.
\]

\begin{rem}
If, in (\ref{tk_kitevocsokkentes}), 
the exponent of $|z_i|$ had not been changed to 2, then the integral of $g$ would have been 
significantly more complicated. Interestingly, similar complication occurs, if one considers simply $|z_i|$ instead of 
$|z_i|^2$.  The rational pair $\frac{1}{4}$ and $\frac{1}{2}$ in the definition of $g$ has
also been a fortunate choice: when working with the numbers $\frac{1}{5}$ and $\frac{1}{2}$ instead, for example, 
\textit{Mathematica} produced so complicated integrals that were practically useless from the viewpoint of further analysis.\quad $\Box$
\end{rem}

\begin{rem}
 An alternative approach to analyzing sum 
(\ref{40-es_tk}) is to estimate $e^{-\frac{h\al}{4}(n-i)}$ above by $1$. However, the resulting integral would not be much simpler in that case either; we would then finally 
use the boundedness of $\al\ln \al$ for $\al\in (0,\al_0]$. Compare the above calculations with their counterparts in the PF case.\quad $\Box$
\end{rem}

Finally, we prove a closeness estimate on $[z_0,0)$ for $\al\le 0$. We begin with a simple observation on monotonicity of the sequence $z_n\equiv z_n(\al)$ (as before, for brevity, the dependence on $h$ is still suppressed).

\begin{lemma}
Suppose that $\al\le 0$ and assumption (\ref{tk_0kicsi}) hold. Then for any $0<h\le h_0$, $-\al_0\le\al\le \beta\le 0$ and $n\in \mathbb{N}$ we have that
\[
0> z_n(\al)\ge z_n(\beta).
\]
\end{lemma}
\textit{Proof.} By definition, we have that $z_0(\al)=z_0(\beta)=z_0$, so suppose that for some  $n$ we already know that $z_n(\al)\ge z_n(\beta)$. Then, by the definition of the sequence $z_n$, and by the facts that the function $z\mapsto \Nfi(h,z,\al)$ is monotone increasing and the function $\al\mapsto \Nfi(h,z,\al)$ is monotone decreasing, we get that
\[
z_{n+1}(\al)=\Nfi(h,z_n(\al),\al)\ge \Nfi(h,z_n(\beta),\al)\ge \Nfi(h,z_n(\beta),\beta)=z_{n+1}(\beta),
\]
which completes the induction.
\quad $\Box$

This means that $0> z_n(\al)\ge z_n(0)$ holds for $\al\le 0$, hence it is enough to give a lower estimate for $z_n(0)$. Such an estimate is presented in \cite[Lemma 2.4.3]{dissertation} with
an easy proof via induction. We just cite the appropriate result without proof.
(In fact, \cite[Lemma 2.4.3]{dissertation} is more general, and that estimate has been formulated in the context of the fold bifurcation. A similar version with proof is
found in \cite[Lemma 3.9]{paezloczisurvey}. See also Lemma \ref{PF414-eslemma}
in the present work.)
\begin{lemma}
Under assumption (\ref{tk_0kicsi}) and for $n\in \mathbb{N}$ we have that
\[
z_n(0)\ge z_0,
\]
and for $n\ge \lfloor \frac{1}{h} \rfloor +1$
\[
z_n(0)\ge -\frac{2}{n h}. 
\]
\end{lemma}
Then we can simply estimate (\ref{tk_outer}) for $\al\le 0$ as follows. Supposing that $n\ge 1$ we get that 
\[
\sup_{[z_n,z_{n+1}]}|id-J^E|\le  c\cdot h^{p+1} |z_0|^3\prod_{j=1}^{n}\left(1+h\al+\frac{7}{4}h\max\left(z_{j},J^E(z_{j})\right)\right)+
\]
\[
c\cdot h^{p+1} \sum_{i=0}^{n-1}|z_i|^3\prod_{j=i+2}^{n}\left(1+h\al+\frac{7}{4}h\max\left(z_{j},J^E(z_{j})\right)\right)\le
\]
\[
c\cdot h^{p+1} |z_0|^3\cdot 1+
c\cdot h^{p}\cdot h \sum_{i=0}^{n}|z_i(0)|^3\cdot 1\le
\]
\[
c\cdot h^{p}\left(h |z_0|^3+h \sum_{i=0}^{\lfloor \frac{1}{h} \rfloor}|z_i(0)|^2+h\sum_{i=\lfloor \frac{1}{h}  \rfloor +1}^{n}|z_i(0)|^2\right),
\]
where, of course, for $n\le \lfloor \frac{1}{h} \rfloor$ the $\sum_{i=\lfloor \frac{1}{h}  \rfloor +1}^{n}$ piece should be omitted. But
\[
h \sum_{i=0}^{\lfloor \frac{1}{h} \rfloor}|z_i(0)|^2 \le (h+1)\cdot z_0^2=2z_0^2,
\]
and
\[
h\sum_{i=\lfloor \frac{1}{h}  \rfloor +1}^{n}|z_i(0)|^2\le h\sum_{i=\lfloor \frac{1}{h} \rfloor +1 }^{n}\frac{4}{i^2 h^2}\le \frac{4}{h}\int_{\frac{1}{h}-1}^{\infty}\frac{1}{i^2}=\frac{4}{1-h}\le 8.
\]
We have thus proved that 
\[
\sup_{[z_0,0)}|id-J^E|\le 12c\cdot h^{p}.
\]

\section{Conjugacy in the Discretized PF Bifurcation}\label{conjugacyPF}

\subsection{Construction of the conjugacy in the PF case}\label{pf_constructionoftheconjugacy}

According to Section \ref{PFmainresultssection}, the normal forms in the PF case can be chosen as
\begin{equation}\label{pf_normalFi}
\NFi(h,x,\al):=(1+h\al)x-h x^3+h x^4\,\eta_4(h,x,\al)
\end{equation}
and
\begin{equation}\label{pf_normalfi}
\Nfi(h,x,\al):=(1+h\al)x-h x^3+h x^4\,{\widetilde{\eta}_4}(h,x,\al),
\end{equation}
where $\eta_4$ and $\widetilde{\eta}_4$ are smooth functions.  
Let $K>0$ denote a uniform bound on $\Big|\frac{\dd^i}{\dx^i}\,\eta(h,\cdot,\al)\Big|$ ($i\in\{0,1,2\}$,\,$\eta\in\{\eta_4,\widetilde{\eta}_4\}$) in a neighborhood of the origin for any
small $h>0$ and $|\al|$, as well as a uniform bound on $\Big|\frac{\dd}{\dd\al}\,\eta(h,x,\cdot)\Big|$ ($\eta\in\{\eta_4,\widetilde{\eta}_4\}$) in a ] of the origin for any
small $h>0$ and $|x|$. 
We also have that there exists a constant $c>0$ such that 
\begin{equation}\label{pf_normalformakkozelsege}
\left| \NFi(h,x,\al)-\Nfi(h,x,\al)\right|\le c\cdot h^{p+1}x^4
\end{equation}
holds for all sufficiently small $h>0$, $|x|\ge 0$ and $|\al|\ge 0$. Throughout Section
\ref{pf_constructionoftheconjugacy}, $c$ will denote this particular positive constant. 
In what follows, suppose that 
\[
0<h\le h_0:=\min\left(\frac{1}{10},8K^2\right),
\]
\begin{equation}\label{pf_0kicsi}
 |x|\le \varepsilon_0:=\min\left(\frac{1}{10},\frac{1}{5K}\right)\ \mathrm{and\ } 
\end{equation} 
\[
|\al| \le \al_0:=\min\left(\frac{1}{288},\frac{1}{72K^2}\right).
\]
However, if the domain of definition of the functions $\eta_4$ and $\widetilde{\eta}_4$ is smaller than $(0,h_0]\times [-\varepsilon_0,\varepsilon_0]\times [-\al_0,\al_0]$ given above, then $h_0$, $\varepsilon_0$ or $\al_0$ should be decreased further.

\begin{figure}[h]
 \centering
 \includegraphics{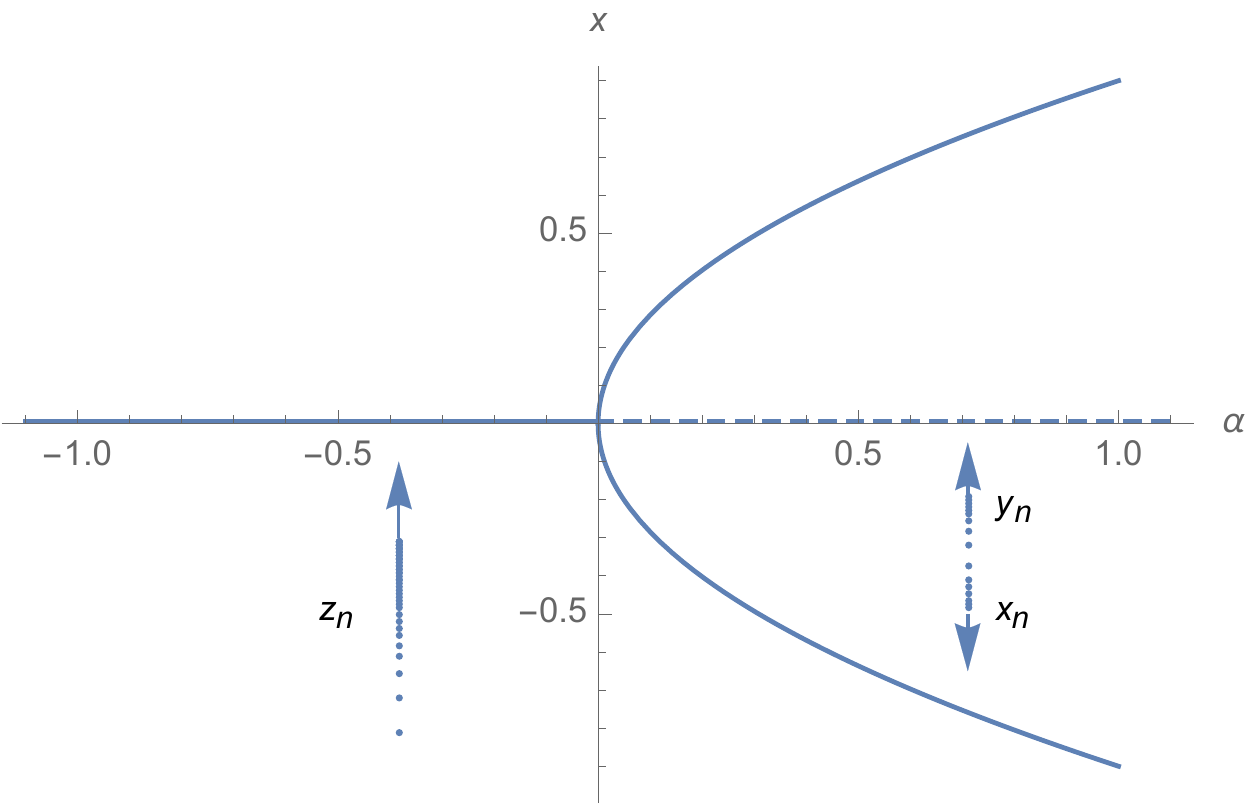}
 \caption{The figure shows the branch of stable (continuous line) and unstable (dashed line) fixed points of $\Nfie$ 
 together with the first few terms of the inner sequences
$x_n(h,\alpha)$ and $y_n(h,\alpha)$, and the outer sequence $z_n(h,\alpha)$, for some 
fixed $h>0$ and $\al$. The arrows indicate the direction of these sequences.}\label{figure4}
\end{figure}

We have that for any fixed $h>0$,
$\omFinul(h,\al)=0$ is an attracting fixed point of the map $\NFi(h,\cdot,\al)$ for $\al\le 0$,
and repelling for $\al>0$, see Figure \ref{figure4}. For any fixed $h>0$ and $\al>0$, this map 
possesses another two attracting fixed points, denoted by $\omFiplusz\equiv\omFiplusz(h,\al)>0$  and $\omFi\equiv\omFi(h,\al)<0$. The three branches of fixed points, $\omFinul(h,\al)$ and $\omega_{\Phi,_{\pm}}(h,\al)$ merge at $\al=0$.
Analogous results hold for the map $\Nfi(h,\cdot,\al)$. Its fixed points are denoted by
$\omfinul \equiv 0$, $\omfi$ and $\omfiplusz$.

 The construction of the homeomorphism $J^E$ is completely analogous to that in the TC case
(the braches of fixed points of the PF and TC bifurcation in the lower half-plane $x\le 0$ look topologically the same), hence it is omitted here. Here, again, we deal only with the $x\le 0$ region---the $x>0$ case is analogous.
The only difference in the construction compared to Section 
\ref{conjugacyTC} is that we set $x_0:=y_0:=-\sqrt{\frac{\al}{8}}$ 
(for $\al>0$) as a starting value in the inner region (but we keep
the same $z_0:=-\varepsilon_0$ starting
value in the outer region). Then conditions (\ref{pf_0kicsi}) imply
that $\Nfi^E$ and $\NFi^E$ (together with their inverses) are strictly  increasing, and $|\al|<\frac{392}{K^2}$ implies  $x_0(\al)>x_1(h,\al)$, while $\sqrt{2\al_0}<|z_0|<\frac{1}{2K}$ implies $z_0<z_1(h,\al)$.  
This means that $x_n$ is strictly decreasing, $y_n$ is strictly increasing (for $\al>0$ and $n\ge 0$), and $\lim_{n\rightarrow \infty}x_n(h,\al)=\omfi$, while $\lim_{n\rightarrow \infty}y_n(h,\al)=\omfinul$.
As for $z_n$, it is strictly increasing, moreover, $\lim_{n\rightarrow \infty}z_n(h,\al)=\omfi$ for $\al>0$, and $\lim_{n\rightarrow \infty}z_n(h,\al)=\omfinul$ for $\al\le 0$.

\subsection{Conjugacy estimates: optimality at the PF fixed points}\label{pf_optimality}

We prove here that the constructed conjugacy $J^E$ is $\cal{O}$$(h^p\,\al)$-close to the
identity at the fixed points $\omfi(h,\al)$, and an explicit example shows
that this estimate is optimal in $h$ and $\al$. This means that the $\cal{O}$$(h^p)$-closeness estimates in Theorem \ref{Conjugacy estimate in the PF bifurcation}
are optimal in $h$.

The first two auxiliary lemmas will frequently be used.

\begin{lemma}\label{pf_Nfi fixpontja}
For every $0<h\le h_0$ and $0< \al \le \al_0$ we have that
\[
\{\omfi,\omFi\} \subset \left(-\sqrt{2\al},-\frac{4}{5}\sqrt{\al}\right ) \subset \left(-\sqrt{2\al},-\sqrt{\frac{3}{5}\al}\right ).
\]
\end{lemma}
\textit{Proof.} By definition, $\omfi<0$ solves
$\al-x^2+x^3\cdot\widetilde{\eta}_4(h,x,\al)=0$, 
that is
\[
\omfi=-\sqrt{\frac{\al}{1-\omfi\cdot\widetilde{\eta}_4(h,\omfi,\al)}}.
\]
But $|x\,\widetilde{\eta}_4|\le \varepsilon_0 K\le \frac{1}{2}$ implies, for example,
$\frac{1}{2} \le 1-\omfi\cdot\widetilde{\eta}_4(h,\omfi,\al)\le \frac{25}{16}$, completing the proof.
In some calculations, the weaker upper bound $-\sqrt{\frac{3}{5}\al}$ will yield a numerically simpler result. The proof for $\omFi$ is similar.
\quad $\Box$

\begin{lemma}\label{pf_NFidervialtbecsles} For any $0<h\le h_0$, $-\varepsilon_0\le x<0$ and $-\al_0\le \al\le \al_0$, we have that
\[
\left(\NFi^E\right)^{\prime}(x)\le 1+h\al-\frac{5}{2}h x^2.
\]
\end{lemma}
\textit{Proof.} The conditions in (\ref{pf_0kicsi}) have been set up to imply this inequality.
\quad $\Box$

\begin{lemma} For any $0<h\le h_0$ and $0<\al\le \al_0$ satisfying (\ref{pf_0kicsi}), we have that
\[
|\,\omfi-\omFi|\le 8 c\cdot h^{p}\,\al.
\]
\end{lemma}
\textit{Proof.} (cf. its TC counterpart) By Lemma \ref{pf_Nfi fixpontja}, (\ref{pf_normalformakkozelsege}) and Lemma \ref{pf_NFidervialtbecsles} we have that
\[|\,\omfi-\omFi|\le 4c\cdot h^{p+1}\al^2+\left(1-\frac{h\al}{2}\right)|\,\omfi-\omFi|,
\]
which yields the desired result.
\quad $\Box$

\begin{rem} The example below shows a situation when the distance of fixed points of normal forms 
satisfying (\ref{pf_normalformakkozelsege}) is bounded from below by
$const \cdot h^p$ (as $h\rightarrow 0$). Since now we are going to deal with cubic polynomials,
we do not attempt to construct their explicit solutions (as we did in the TC case having quadratic polynomials), but we give a more general, yet simpler argument.
Set 
$\NFi(h,x,\al):=(1+h\al)x-h x^3$ and $\Nfi(h,x,\al):=(1+h\al)x-h x^3+h^{p+1}x^4$. Then these maps satisfy
(\ref{pf_normalformakkozelsege}) in a neighborhood of the origin, moreover, $\omFi=-\sqrt{\al}$.
As for $\omfi$, we see that $\omfi=-\sqrt{\frac{\al}{1-h^p\,\omfi}}$. Then by Lemma \ref{pf_Nfi fixpontja} we get that 
\[
\omfi\in \left(-\sqrt{\frac{\al}{1+h^p\,\frac{4}{5}\sqrt{\al}}},-\sqrt{\frac{\al}{1+h^p\,\sqrt{2\al}}}\right)
\]
so $\omfi>-\sqrt{\al}=\omFi$, yielding
\[
|\,\omfi-\omFi|\ge \left| -\sqrt{\al}+ \sqrt{\frac{\al}{1+h^p\,\frac{4}{5}\sqrt{\al}}}\right|=
\sqrt{\al} \left| 1- \frac{1}{\sqrt{1+t}}\right|
\]
with $t:=h^p\,\frac{4}{5}\sqrt{\al}$. Then, by (\ref{pf_0kicsi}), $t\in (0,1)$. But
for any such $t$ 
\[
\left| 1- \frac{1}{\sqrt{1+t}}\right|=\frac{t}{\sqrt{1+t}(1+\sqrt{1+t})}\ge \frac{t}{4},
\]
hence
\[ 
|\,\omfi-\omFi|\ge \frac{1}{5} h^p\, \al.
\quad \Box
\]
\end{rem}

\subsection{Conjugacy estimates in the inner PF region}\label{pf_inner region section}\label{pf_innerregion}

In this section, for any fixed $0<h\le h_0$ and $0<\al\le \al_0$, the closeness estimate is proved in $(\omfi,0)$ in a similar way as in the TC case, hence most intermediate steps
and inequalities are omitted or only sketched (however, the key Lemma \ref{pf_xknagysagrend}
 is carefully examined). First, the proof is presented in $(\omfi,x_0]$ (Lemma \ref{pflemma4.8}),
then in $(x_0,0)$ (Lemma \ref{pflemma4.10}).

Now we have that
\[
\sup_{[x_1,x_{0}]}|id-J^E|= \frac{c}{64}h^{p+1}\al^2,
\]
while for $n\ge 1$
\[
\sup_{[x_{n+1},x_n]}|id-J^E|\le  \sup_{[x_n,x_{n-1}]}\left|\Nfi^E-\NFi^E\right|+ \sup_{x\in [x_n,x_{n-1}]} \left(\left(\sup_{[\{x,J^E(x)\}]}
\left(\NFi^E\right)^{\prime}\right)|x-J^E(x)|\right)\le  \]
\[
c\cdot h^{p+1}  x_n^4+ \left(1+h\al-\frac{5}{2}h\max\left(x_{n-1},J^E(x_{n-1})\right)^2\right) \sup_{[x_n,x_{n-1}]}|id-J^E|,
\]
where we have used Lemma \ref{pf_NFidervialtbecsles}, the fact that the functions $id$ and $J^E$ are increasing, and inequality
\[
\sup_{[\{x,J^E(x)\}]}\left(\NFi^E\right)^{\prime}\le \sup_{[\{x,J^E(x)\}]}\left(1+h\al-\frac{5}{2}h\cdot id^2\right)\le 1+h\al-\frac{5}{2}h\max\left(x,J^E(x)\right)^2.
\]
In order to prove that the conjugacy $J^E$ is $\cal{O}$$(h^p)$-close to the identity
on $(\omfi,x_0]$ for any $h\in (0,h_0]$ and $\al \in (0,\al_0]$, we will show that
\begin{equation}\label{pf_innerbecslese}
\sup_{h\in (0,h_0]}\sup_{\al\in (0,\al_0]}\sup_{n\in \mathbb{N}} h \sum_{i=0}^{n} 
x_i^4\prod_{j=i}^{n-1}\left(1+h\al-\frac{5}{2}h\max\left(x_{j},J^E(x_{j})\right)^2\right)\le const
\end{equation}
holds with a suitable $const\ge 0$ (where $\prod_{j=n}^{n-1}$ is understood to be $1$).

First an explicit upper estimate of the sequence $\max\left(x_{n},J^E(x_{n})\right)$ is given. 

\begin{lemma}\label{pf_xknagysagrend}
For $n\ge 0$, set \[a_n(h,\al):=-\frac{4}{5}\sqrt{\al}\cdot\frac{(1+h\al)^{n}}{\sqrt{5+(1+h\al)^{2n}}},\] then we have that
$x_n\in (\omfi,a_n)$ and $J^E(x_n)\in (\omFi,a_n)$.
\end{lemma}
\textit{Proof.} Due to assumptions (\ref{pf_0kicsi}), $\max\left(\omfi,\omFi\right)<a_n$ for $n\ge 0$, so the intervals in the lemma are non-degenerate. We proceed by induction. We see that $a_0>x_0\equiv J^E(x_0)\equiv -\sqrt{\frac{\al}{8}}$ is always satisfied.  
So suppose that the statement is true for some $n\ge 0$. Condition $|x|<\frac{1}{3K}$
implies $\Nfie(x)<(1+h\al)x-\frac{4}{3}h x^3$, moreover, by monotonicity of $\Nfi^E$ we get that
\begin{equation}\label{pf_subtle}
x_{n+1}=\Nfie(x_n)<\Nfie(a_n)<(1+h\al)a_n-\frac{4}{3}h a_n^3.
\end{equation}
Therefore it is enough to prove that 
\begin{equation}\label{pf_anpozitive}
a_{n+1}-(1+h\al)a_n+\frac{4}{3}h a_n^3>0. 
\end{equation}
For brevity, we set $\lambda:=h\al>0$. Then (\ref{pf_anpozitive}) is equivalent to
\[
\frac{4}{375}\left( -A + B - C \right)\sqrt{\alpha} \left( 1 + \lambda  \right) ^n > 0,
\] where $A:=\frac{64\,\lambda \,{\left( 1 + \lambda  \right) }^{2\,n}}
  {{\left( 5 + {\left( 1 + \lambda  \right) }^{2\,n} \right) }^
    {\frac{3}{2}}}$, $B:=\frac{75\,\left( 1 + \lambda  \right) }
  {{\sqrt{5 + {\left( 1 + \lambda  \right) }^{2\,n}}}}$ and $C:=\frac{75\,\left( 1 + \lambda  \right) }
  {{\sqrt{5 + {\left( 1 + \lambda  \right) }^
        {2 + 2n }}}}$. We will show that \[ -A + B - C >0.\]

First put $B-C$ over a common denominator. Then, to eliminate square roots from its numerator, multiply it by
$\frac{{\sqrt{5 + {\left( 1 + \lambda  \right) }^{2\,n}}} + 
  {\sqrt{5 + {\left( 1 + \lambda  \right) }^
       {2 + 2\,n}}}}{{\sqrt{5 + {\left( 1 + \lambda  \right) }^{2\,n}}} + 
  {\sqrt{5 + {\left( 1 + \lambda  \right) }^
       {2 + 2\,n}}}}$. After these manipulations, the product $\lambda(1+\lambda)^{2n}>0$ can be factored out from all three terms. Hence
$-A+B-C>0$ becomes 
\[
\frac{-64}{{\left( 5 + 
        {\left( 1 + \lambda  \right) }^{2\,n} \right) }
      ^{\frac{3}{2}}} + 
\]
\[
  \frac{75\,\left( 1 + \lambda  \right) \,
     \left( 2 + \lambda  \right) }{{\sqrt{5 + 
         {\left( 1 + \lambda  \right) }^{2\,n}}}\,
     {\sqrt{5 + {\left( 1 + \lambda  \right) }^
          {2 + 2\,n}}}\,
     \left( {\sqrt{5 + 
           {\left( 1 + \lambda  \right) }^{2\,n}}} + 
       {\sqrt{5 + {\left( 1 + \lambda  \right) }^
            {2 + 2\,n}}} \right) }>0.
\]
The left hand side of the above inequality is decreased, if the denominator of the second term
is replaced by $2\,{\sqrt{5 + {\left( 1 + \lambda  \right) }^{2\,n}}}\,
  \left( 5 + {\left( 1 + \lambda  \right) }^
     {2 + 2\,n} \right)$. Then we can multiply the expression by ${\sqrt{5 + {\left( 1 + \lambda  \right) }^{2\,n}}}$ and all square roots are got rid of: we are to verify
\[\frac{-64}{5 + {\left( 1 + \lambda  \right) }^
      {2\,n}} + \frac{75\,
     \left( 1 + \lambda  \right) \,
     \left( 2 + \lambda  \right) }{2\,
     \left( 5 + {\left( 1 + \lambda  \right) }^
        {2 + 2\,n} \right) }>0.
\]
By condition $(\ref{pf_0kicsi})$, $(1+\lambda)^2<\frac{75}{64}$, so it is enough to show that 
\[
\frac{-64}{5 + {\left( 1 + \lambda  \right) }^
       {2\,n}} + \frac{75\,
      \left( 1 + \lambda  \right) \,
      \left( 2 + \lambda  \right) }{2\,
      \left( 5 + \frac{75\,
           {\left( 1 + \lambda  \right) }^{2\,n}}{64}
        \right) } > 0.
\]
However, the left hand side above can be factored to yield
\[
\frac{32\,\left( 22 + 225\,\lambda  + 
      75\,{\lambda }^2 + 
      45\,\lambda \,{\left( 1 + \lambda  \right) }^
        {2\,n} + 15\,{\lambda }^2\,
       {\left( 1 + \lambda  \right) }^{2\,n} \right) }
    {\left( 5 + {\left( 1 + \lambda  \right) }^
       {2\,n} \right) \,
    \left( 64 + 15\,{\left( 1 + \lambda  \right) }^
        {2\,n} \right) },
\]
which is clearly positive. 

The proof for the sequence $J^E(x_{n})$ is the completely similar: by construction of $J$, the beginning of (\ref{pf_subtle}) should (and can) be replaced by $J^E(x_{n+1})=\NFi^E(J^E(x_n))<\NFi(a_n)$, but then every step is unchanged.
\quad $\Box$

\begin{rem} Attempts to approximate 
subexpressions of the form $(a+b t)^{\gamma}$ with their series expansions (up to third order)  turned out to be insufficient to complete the proof. To find the above ``purely algebraical" manipulations, \textit{Mathematica} has been extensively used. The definition of $a_n$ is again based on the beautiful parameterized model function (\ref{thulsmodelfunction}).\quad $\Box$
\end{rem}

The sum $\sum_{i=0}^n$ in (\ref{pf_innerbecslese}) is split into two at $\lceil \frac{6} {h\al}\rceil$. This choice is motivated by the following lemma.

\begin{lemma}\label{pf_inner split}
Suppose that $n>\lceil \frac{6}{h\al}\rceil$. Then
\[
\max\left(x_n,J^E(x_n)\right)<-\sqrt{\frac{3}{5}\al},
\]
hence 
\[
1+h\al-\frac{5}{2}h\max\left(x_n,J^E(x_n)\right)^2<1-\frac{h\al}{2}
\]
holds for $n>\lceil \frac{6}{h\al}\rceil$.
\end{lemma}
\textit{Proof.} By virtue of Lemma \ref{pf_xknagysagrend} it is enough to show that
$n>\lceil \frac{6}{h\al}\rceil$ implies $a_n<-\sqrt{\frac{3}{5}\al}$. But this latter  inequality is equivalent to $(1+h\al)^n >\sqrt{75}$. However, $e^3>\sqrt{75}$, so the  corresponding proof given in the TC case suffices here, too. 
\quad $\Box$

Let us turn directly to (\ref{pf_innerbecslese}) now and fix $h\in (0,h_0]$, $\al\in (0 ,\al_0]$ and $n\in \mathbb{N}^+$ arbitrarily. (If
$n\le \lceil \frac{6}{h\al}\rceil$, then the sums $\sum_{i=\lceil \frac{6}{h\al}\rceil+1}^n$  below are absent and the proof is simpler.) Using $\omfi<x_i<0$,  $|x_i|\le \sqrt{2\al}$ (by  Lemma \ref{pf_Nfi fixpontja}), and $\max\left(x_j,J^E(x_j)\right)\le x_0\equiv J^E(x_0) \equiv -\sqrt{\frac{\al}{8}}$ (by monotonicity), moreover, Lemma \ref{pf_inner split},  assumption $h\al<1$ from (\ref{pf_0kicsi}) and inequality $(1+\frac{1}{A})^A\le e$ (for  $A\ge 1$), we see that
\[
h \sum_{i=0}^{n} 
x_i^4\prod_{j=i}^{n-1}\left(1+h\al-\frac{5}{2}h\max\left(x_{j},J^E(x_{j})\right)^2\right)\le 
\]
\[
4h\al^2\sum_{i=0}^{\lceil \frac{6}{h\al}\rceil} 
\prod_{j=1}^{\lceil \frac{6}{h\al}\rceil-1}\left(1+h\al-\frac{5}{2}\cdot\frac{h\al}{8} \right) 
+4h\al^2\sum_{i=\lceil \frac{6}{h\al}\rceil+1}^{n}\prod_{j=i}^{n-1}\left(1-\frac{h\al}{2} \right)\le 
\]
\[
4h\al^2\left(1+\frac{11}{16}h\al\right)^{\frac{6}{h\al}}\left(\Big{\lceil}{\frac{6}{h\al}} \Big{\rceil}+1\right)+
4h\al^2\sum_{i=\lceil \frac{6}{h\al}\rceil+1}^{n}\left(1-\frac{h\al}{2}\right)^{n-i}\le
\]
\[
4h\al^2\left(1+\frac{11}{16}h\al\right)^{\frac{16}{11h\al}\cdot\frac{11h\al}{16}\cdot\frac{ 6}{h\al}}\left(\frac{6+2h\al}{h\al}\right)+
4h\al^2\sum_{i=0}^{\infty}\left(1-\frac{h\al}{2}\right)^{i}\le
\]
\[
4h\al^2\cdot e^{\frac{66}{16}}\cdot \frac{8}{h\al}+
4h\al^2\cdot\frac{2}{h\al}\le 1988\,\al.
\]

\noindent Therefore, $\sup_{[x_{n+1},x_{n}]}|id-J^E|\le 1988c\cdot h^p \al$  for any $h\in  (0,h_0]$, $\al\in (0,\al_0]$ and $n\ge 1$, and, as we have seen, $\sup_{[x_{1},x_{0}]} |id-J^E|= \frac{c}{64}h^{p+1}\al^2$, which yield the following lemma.
\begin{lemma}\label{pflemma4.8}
Under assumption (\ref{pf_0kicsi})
\[
\sup_{(\omfi,x_0]}|id-J^E|\le 1988c\cdot h^p \al.
\]
\end{lemma}

In the rest of the section, the closeness estimate is proved in the interval $(y_0,\omfinul)$. Recall that $y_0=x_0=J^E(x_0)\equiv -\sqrt{\frac{\al}{8}}$ and $\omfinul=\omFinul\equiv 0$.

Suppose first that $n\ge 1$ (the case $n=0$ will be examined later). Then we proceed
exactly as in the TC case, so we will only list the differences. We get that
\[
\sup_{[y_n,y_{n+1}]}|id-J^E|\le 
\] 
\[
\left[ \sup_{x\in [y_n,y_{n+1}]}\sup_{[\{\NFi^E(x),J^E\circ \Nfi^E(x)\}]}\left(\left(\NFi^E\right)^{[-1]}\right)^{\prime}\right]
\left[c\cdot h^{p+1} y_n^4+ \sup_{[y_{n-1},y_{n}]} |id-J^E|\right]. 
\]

\noindent The following lemma gives an upper bound on the first term above (and shows 
a motivation for the choice of $x_0=-\sqrt{\frac{\al}{8}}$).

\begin{lemma}
Suppose that $n\ge 1$, then under assumption (\ref{pf_0kicsi}) we have that
\[
 \sup_{x\in [y_n,y_{n+1}]}\sup_{[\{\NFi^E(x),J^E\circ \Nfi^E(x)\}]}\left(\left(\NFi^E\right)^{[-1]}\right)^{\prime}\le 1-\frac{h\al}{4}.
\]
\end{lemma}
\textit{Proof.} As in the TC case, we have that
\[
\sup_{x\in [y_n,y_{n+1}]}\sup_{[\{\NFi^E(x),J^E\circ \Nfi^E(x)\}]}\left(\left(\NFi^E\right)^{[-1]}\right)^{\prime}\le
\sup_{\left(-\sqrt{\frac{\al}{8}},0\right)}\frac{1}{\left(\NFi^E\right)^{\prime}}\le \ldots
\]
But assumption (\ref{pf_0kicsi}) together with $x<0$ imply that $\left(\NFi^E\right)^{\prime}(x)\ge
1+h\al-4hx^2\ge 0$. So
\[
\ldots \le \sup_{x\in \left(-\sqrt{\frac{\al}{8}},0\right)}\frac{1}{1+h\al-4h x^2}\le \frac{1}{1+h\al-4h \left(-\sqrt{\frac{\al}{8}}\right)^2}=\frac{1}{1+\frac{1}{2}h\al}\le 1-\frac{h\al}{4}.
\]
\quad $\Box$

We have thus proved (using $|y_n|\le \sqrt{\frac{\al}{8}}$ also) that for $n\ge 1$
\[
\sup_{[y_n,y_{n+1}]}|id-J^E|\le \left(1-\frac{h\al}{4}\right)\left[\frac{c}{64}\cdot h^{p+1} \al^2+ \sup_{[y_{n-1},y_{n}]} |id-J^E|\right].
\]
For $n=0$, similarly as in the TC case, we get that
\[
\sup_{[y_0,y_{1}]}|id-J^E|\le 2\cdot \frac{c}{64}h^{p+1}\al^2,
\]
and for $n\ge 1$ that
\[
\sup_{[y_n,y_{n+1}]}|id-J^E|\le \left(1-\frac{h\al}{4} \right)^n \sup_{[y_0,y_{1}]}|id-J^E| + \frac{c}{64}h^{p+1}\al^2\sum_{i=1}^{n} \left(1-\frac{h\al}{4} \right)^i\le
\]
\[
1\cdot 2\cdot \frac{c}{64}h^{p+1}\al^2+\frac{c}{64}h^{p+1}\al^2 \cdot \frac{4}{h\al}\le \frac{c}{8}h^{p}\al, 
\]
using $h\al\le 1$ by  (\ref{pf_0kicsi}). Since the same upper estimate is valid for $n=0$, too, we have proved the following lemma.
\begin{lemma}\label{pflemma4.10}
Under assumption (\ref{pf_0kicsi})
\[
\sup_{(x_0,0)}|id-J^E|\le \frac{c}{8}h^{p}\al.
\]
\end{lemma}

\subsection{Conjugacy estimates in the outer PF region}\label{pf_outer_label}\label{pf_outerregion}

In this section, we first prove an ${\cal{O}}(h^p)$-closeness estimate in the interval $[z_0,\omfi)$ for $\al>0$, then the closeness is proved on $[z_0,\omFinul)\equiv [z_0,0)$ for $\al\le 0$. The
key lemma of this section is Lemma \ref{pf_zknagysagrend}.

We are already familiar with the inequalities below (cf. the TC case). For $n\ge 1$ we have that
\[
\sup_{[z_n,z_{n+1}]}|id-J^E|\le c\cdot h^{p+1} z_0^4\prod_{j=1}^{n}\left(1+h\al-\frac{5}{2}h\max\left(z_{j},J^E(z_{j})\right)^2\right)+
\]
\begin{equation}\label{pf_outer}
c\cdot h^{p+1} \sum_{i=0}^{n-1}z_i^4\prod_{j=i+2}^{n}\left(1+h\al-\frac{5}{2}h\max\left(z_{j},J^E(z_{j})\right)^2\right),
\end{equation}
where $\prod_{j={n+1}}^{n}$ is, as always, 1, and 
\[
\sup_{[z_0,z_{1}]}|id-J^E|\le c\cdot h^{p+1} z_0^4.
\]

 The lemma gives a lower estimate of the sequence $z_n$ for $\al>0$.

\begin{lemma}\label{pf_zknagysagrend}
For $n\ge 0$, set \[b_n(h,\al):=-2\sqrt{\al}\cdot\frac{(1+h\al)^{n}}{\sqrt{\al-1+(1+h\al)^{2n}}},\] then 
$b_n\le \min\left(z_n,J^E(z_n)\right)$.
\end{lemma}
\textit{Proof.} We prove by induction: $b_0=-2< z_0=J^E(z_0)$ holds due to assumption (\ref{pf_0kicsi}). So suppose that the statement is true for some $n\ge 0$. 
We have that $\Nfie(x)\ge (1+h\al)x-\frac{3}{4}h x^3$ (by $x<0$ and $|x|\le \varepsilon_0<\frac{1}{4K}$), moreover, that $(1+h\al)id-\frac{3}{4}h\,id^3$ is monotone increasing (which is true, for example, for $|x|\le \frac{2}{3\sqrt{h}}$; but $|b_n|\le 2$ is easily seen, and due to 
$h\le\frac{1}{10}$ we get $|b_n|\le \frac{2}{3\sqrt{h}}$ also), 
so we obtain that
\begin{equation}\label{pf_subtle2}
z_{n+1}=\Nfie(z_n)\ge (1+h\al)z_n-\frac{3}{4}h z_n^3\ge (1+h\al)b_n-\frac{3}{4}h b_n^3. 
\end{equation}
Therefore it is sufficient to show that 
\[
(1+h\al)b_n-\frac{3}{4}h b_n^3\ge b_{n+1},
\]
being the same as
\[
0\le 2\,{\sqrt{\alpha }}\,{\left( 1 + \lambda  \right) }^n\,
   \left(  \widetilde{A}-\widetilde{B}+\widetilde{C} \right),\]
with $\widetilde{A}:=\frac{3\,\lambda \,{\left( 1 + \lambda  \right) }^{2\,n}}
      {{\left( -1 + \alpha  + {\left( 1 + \lambda  \right) }^{2\,n} \right) }^{\frac{3}{2}}}$, $\widetilde{B}:=\frac{1 + \lambda }{{\sqrt{-1 + \alpha  + {\left( 1 + \lambda  \right) }^{2\,n}}}}$,  $\widetilde{C}:=\frac{1 + \lambda }{{\sqrt{-1 + \alpha  + {\left( 1 + \lambda  \right) }^{2 + 2\,n}}}}$ and $\lambda:=h\al >0$. 

Now proceeding just as in Lemma \ref{pf_xknagysagrend}, we first get
\[
0\le\frac{3}{{\left( -1 + \alpha  + {\left( 1 + \lambda  \right) }^{2\,n} \right) }^
     {\frac{3}{2}}} -
\]
\[\frac{\widetilde{B}\,\left( 2 + \lambda  \right) }
   {{\sqrt{-1 + \alpha  + {\left( 1 + \lambda  \right) }^{2 + 2\,n}}}\,
     \left( {\sqrt{-1 + \alpha  + {\left( 1 + \lambda  \right) }^{2\,n}}} + 
       {\sqrt{-1 + \alpha  + {\left( 1 + \lambda  \right) }^{2 + 2\,n}}} \right) }
\]
to verify. Then multiply the inequality by ${\sqrt{-1 + \alpha  + {\left( 1 + \lambda  \right) }^{2\,n}}}$ to get
\[
0\le \frac{3}{-1 + \alpha  + {\left( 1 + \lambda  \right) }^{2\,n}} -
\]
\[
  \frac{\left( 1 + \lambda  \right) \,\left( 2 + \lambda  \right) }
   {{\sqrt{-1 + \alpha  + {\left( 1 + \lambda  \right) }^{2 + 2\,n}}}\,
     \left( {\sqrt{-1 + \alpha  + {\left( 1 + \lambda  \right) }^{2\,n}}} + 
       {\sqrt{-1 + \alpha  + {\left( 1 + \lambda  \right) }^{2 + 2\,n}}} \right) }.
\]
A sufficient condition for this is
\[
0\le \frac{3}{-1 + \alpha  + {\left( 1 + \lambda  \right) }^{2\,n}} - 
  \frac{\left( 1 + \lambda  \right) \,\left( 2 + \lambda  \right) }
   {2\,\left( -1 + \alpha  + {\left( 1 + \lambda  \right) }^{2\,n} \right) },
\]
but the right hand side is equal to
\[\frac{\left( 1 - \lambda  \right) \,\left( 4 + \lambda  \right) }
  {2\,\left( -1 + \alpha  + {\left( 1 + \lambda  \right) }^{2\,n} \right) },\]
which is positive for $0<\lambda\equiv h\al<1$.

When $\NFi$ and $J^E(z_n)$ are written instead of $\Nfi$ and $z_n$, respectively, the considerations above remain valid, implying $b_n\le J^E(z_n)$. 
\quad $\Box$

\begin{rem} Notice the subtle difference between the chain of inequalities
(\ref{pf_subtle}) and (\ref{pf_subtle2}). Unlike $\Nfie(a_n)$, quantity $\Nfie(b_n)$ is not necessarily defined, since $b_n$ may lie
outside the domain of definition of $\Nfie$.\quad $\Box$
\end{rem}

Now, since $z_j<\omfi<0$ and $J^E(z_j)<\omFi<0$, by Lemma \ref{pf_Nfi fixpontja} we get that the 
right-hand side of (\ref{pf_outer}) is at most 
\[
c\cdot h^{p+1} z_0^4\prod_{j=1}^{n}\left(1-\frac{h\al}{2}\right)+ c\cdot h^{p+1} \sum_{i=0}^{n-1}z_i^4\prod_{j=i+2}^{n}\left(1-\frac{h\al}{2}\right)\le
\]
\[
c\cdot h^{p+1} z_0^4 + c\cdot h^{p+1} \sum_{i=0}^{n-1}z_i^4\left(1-\frac{h\al}{2}\right)^{n-1-i}.
\]
We will show that $h\sum_{i=0}^{n}z_i^4\left(1-\frac{h\al}{2}\right)^{n-i}$ is uniformly bounded for any $n\ge 0$,
$0<h\le h_0$ and $0<\al\le \al_0$.

For $n\ge\lceil\frac{1}{h\al}\rceil$ and $i\ge n$, it is easy to see that $(1+h\al)^i\ge (1+h\al)^{\frac{1}{h\al}}\ge 1+\frac{1}{h\al}\cdot h\al=2$ implies $\frac{(1+h\al)^{i}}{\sqrt{\al-1+(1+h\al)^{2i}}}\le 2$, hence 
by Lemma \ref{pf_zknagysagrend} we have 
\[
h\sum_{i=\lceil\frac{1}{h\al}\rceil}^{n}z_i^4\left(1-\frac{h\al}{2}\right)^{n-i}\le
h\sum_{i=\lceil\frac{1}{h\al}\rceil}^{n}b_i^4\left(1-\frac{h\al}{2}\right)^{n-i}\le 
\]
\[
16h\al^2\sum_{i=\lceil\frac{1}{h\al}\rceil}^{n}\left(\frac{(1+h\al)^{i}}{\sqrt{\al-1+(1+h\al)^{2i}}}\right)^4\left(1-\frac{h\al}{2}\right)^{n-i}\le 256 h \al^2 \sum_{i=0}^{\infty}\left(1-\frac{h\al}{2}\right)^{i} = 512\al.
\]
On the other hand,  due to inequalities $e^{\frac{x}{2}}\le 1+x$ ($x\in [0,1]$) and $1+x\le e^x$ ($x\in \mathbb{R}$), for $n<\lceil\frac{1}{h\al}\rceil$ we have that  
\[
h\sum_{i=0}^{n}z_i^4\left(1-\frac{h\al}{2}\right)^{n-i}\le 
16 h \al^2\sum_{i=0}^{n}\left(\frac{(1+h\al)^{i}}{\sqrt{\al-1+(1+h\al)^{2i}}}\right)^4\left(1-\frac{h\al}{2}\right)^{n-i}\le \]
\[ 
16 h \al^2\sum_{i=0}^{n}\frac{e^{4h\al i}}{\left(\al-1+(1+h\al)^{2i}\right)^2}\cdot 1 \le 
874\, h \al^2\sum_{i=0}^{n}\frac{1}{\left(\al-1+e^{h\al i}\right)^2}. 
\]
Set $g_{h,\al}(x)\equiv g(x):=\frac{h \al^2}{\left(\al-1+e^{h\al x}\right)^2}$ for $x\in [0,\infty)$. Notice that $g$ is bounded at $x=0$, and strictly decreasing on  $[0,\infty)$, because
\[
g^\prime(x)=\frac{-2h^2{\alpha }^3e^{h\al x}}{{\left(\alpha -1 + e^{h\alpha x }   \right) }^    3}<0.
\]
So, since $0<\al<1$, we see that
\[
874\, h \al^2\sum_{i=0}^{n}\frac{1}{\left(\al-1+e^{h\al i}\right)^2}=874h+874\sum_{i=1}^{n}g_{h,\al}(i)\le
874h + 874 \int_{0}^{\frac{1}{h\al}}g_{h,\al}(x)\dx=
\]
\[
874h + 874 \left[
\frac{h\al^2 x}{{\left( \al-1   \right) }^2} + 
  \frac{\alpha }{\left( \al-1  \right) \,
     \left( \al-1 + e^{h\al x}   \right) } - 
  \frac{\alpha \,\ln (\al-1 + e^{h\al x } )}{{\left( \al-1   \right) }^2}
\right]_{x=0}^{\frac{1}{h\al}}=
\]
\[
874h+874\left(
\frac{e-1 + {\alpha }^2 + \alpha \,\left( e-1 + \alpha  \right) \,
     \left( \ln \alpha  - \ln (e-1  + \alpha ) \right) }{{\left( \al-1  \right) }^
     2\,\left( e-1 + \alpha  \right) }
\right)\le
\]
\[
874h+874\cdot
\frac{e-1 + {\alpha }^2   }{{\left( \al-1  \right) }^
     2\,\left( e-1 + \alpha  \right) }\le 
874h+\frac{874 }{{\left( \al-1  \right) }^2 }<3584,
\]
since $h\le \frac{1}{10}$ and $\al\le \frac{1}{2}$, so $\al^2\le \al$.

Now combining all the estimates in the section, under assumption (\ref{pf_0kicsi}) and for $\al>0$
we get that 
\[
\sup_{[z_0,\omfi)}|id-J^E|=  \sup_{n\in \mathbb{N}} \sup_{[z_{n},z_{n+1}]}|id-J^E|\le
\]
\[
\sup_{n\in \mathbb{N}} \max\left(c\cdot h^{p+1} z_0^4,\,c\cdot h^{p+1} z_0^4 + c\cdot h^{p+1} \sum_{i=0}^{n}z_i^4\left(1-\frac{h\al}{2}\right)^{n-i}\right)\le
\]
\[
c\cdot h^{p+1} z_0^4 + c\cdot h^{p}\cdot (3584+512\al)\le 3841 c\cdot h^p.
\]

Finally, a closeness estimate on $[z_0,0)$ for $\al\le 0$ is proved. The proof of the next lemma is identical to its counterpart in the TC case.

\begin{lemma}
Suppose that $\al\le 0$ and assumption (\ref{pf_0kicsi}) hold. Then for any $0<h\le h_0$, $-\al_0\le\al\le \beta\le 0$ and $n\in \mathbb{N}$ we have that
\[
0> z_n(\al)\ge z_n(\beta).
\]
\end{lemma}

 So $0> z_n(\al)\ge z_n(0)$ holds for $\al\le 0$ and it is enough to give a lower estimate of $z_n(0)$.
\begin{lemma}\label{PF414-eslemma}
Under assumption (\ref{pf_0kicsi}), we have for $n\in \mathbb{N}$ that
\[
z_n(0)\ge z_0\ge -\frac{1}{2K}
\]
and for $n\ge \lfloor \frac{16K^2}{h} \rfloor$
\[
z_n(0)\ge -\frac{2}{\sqrt{n h}}.
\]
\end{lemma}
\textit{Proof.} Monotonicity of the sequence $z_n(0)$ and (\ref{pf_0kicsi}) imply the first part. The nontrivial second inequality is proved by induction.  Since 
$z_{\lfloor \frac{16K^2}{h} \rfloor}(0)\ge -\frac{1}{2K}=-\frac{2}{\sqrt{16K^2}}\ge 
-\frac{2}{\sqrt{h \lfloor \frac{16K^2}{h} \rfloor}}$, 
the induction can be started. So suppose that $n\ge \lfloor \frac{16K^2}{h} \rfloor $.
The function $z\mapsto \Nfi(h,z,0)$ is increasing, so, by the induction hypothesis we have that 
$z_{n+1}(0)=\Nfi(h,z_n(0),0)\ge \Nfi\left(h,-\frac{2}{\sqrt{n h}},0\right)$. It is enough to show that  $\Nfi\left(h,-\frac{2}{\sqrt{n h}},0\right)\ge -\frac{2}{\sqrt{(n+1) h}}$. Multiplying this by $\sqrt{h}$ and
rearranging, it suffices to prove 
\[
\frac{4}{n\sqrt{n}}+\frac{8 {\widetilde{\eta}_4}}{n^2\sqrt{h}}\ge \frac{1}{\sqrt{n}}-\frac{1}{\sqrt{n+1}},
\]
where is $\widetilde{\eta}_4$ from (\ref{pf_normalfi}). But $\frac{1}{2n\sqrt{n}}\ge \frac{1}{\sqrt{n}}-\frac{1}{\sqrt{n+1}}$, so it is enough to verify $
4+\frac{8 {\widetilde{\eta}_4}}{\sqrt{n h}}\ge \frac{1}{2}$: it is easy to see that conditions 
$n\ge \lfloor \frac{16K^2}{h} \rfloor \ge \frac{16K^2}{h}-1$ and $h\le 8K^2$ (from (\ref{pf_0kicsi})) together with 
the definition of $K$ imply $\left |\frac{8 {\widetilde{\eta}_4}}{\sqrt{n h}}\right |\le 3$. 
\quad $\Box$

Then we can simply estimate (\ref{pf_outer}) for $\al\le 0$ as follows. Supposing that $n\ge 1$ we get that 
\[
\sup_{[z_n,z_{n+1}]}|id-J^E|\le  c\cdot h^{p+1} z_0^4\prod_{j=1}^{n}\left(1+h\al-\frac{5}{2}h\max \left(z_{j},J^E(z_{j})\right)^2\right)+
\]
\[
c\cdot h^{p+1} \sum_{i=0}^{n-1}z_i^4\prod_{j=i+2}^{n}\left(1+h\al-\frac{5}{2}h\max\left(z_{j},J^E(z_{j})\right)^2\right)\le
\]
\[
c\cdot h^{p+1} z_0^4\cdot 1+
c\cdot h^{p}\cdot h \sum_{i=0}^{n}z_i(0)^4\cdot 1\le
\]
\[
c\cdot h^{p}\left(h z_0^4+h \sum_{i=0}^{\lfloor \frac{16K^2}{h} \rfloor}z_i(0)^4+h\sum_{i=\lfloor \frac{16K^2}{h}  \rfloor +1}^{n}z_i(0)^4\right),
\]
where, of course, for $n \le \lfloor \frac{16K^2}{h} \rfloor$, the sum above $\sum_{i=\lfloor \frac{16K^2}{h}  \rfloor +1}^{n}$ is not present. But
\[
h \sum_{i=0}^{\lfloor \frac{16K^2}{h} \rfloor}z_i(0)^4 \le h\cdot
\left(\frac{16K^2}{h}+ 1 \right)\left(-\frac{1}{2K}\right)^4\le \frac{2}{K^2}
\]
by $h\le 8K^2$, and
\[
h\sum_{i=\lfloor \frac{16K^2}{h}  \rfloor +1}^{n}z_i(0)^4\le h\sum_{i=\lfloor \frac{16K^2}{h} \rfloor +1 }^{n}\frac{16}{i^2 h^2}\le \frac{16}{h}\int_{\frac{16K^2}{h}}^{\infty}\frac{1}{i^2}=\frac{1}{K^2}.
\]
We have thus proved that 
\[
\sup_{[z_0,0)}|id-J^E|\le c\cdot h^{p+1} z_0^4+c\cdot h^{p}\left(h z_0^4+\frac{3}{K^2}\right)\le c \left(2+\frac{3}{K^2}\right)h^{p}.
\]

\section{A numerical example and some open questions}\label{Section5}

In this section we present a TC model example to illustrate the transformations 
in Lemma \ref{Phinormalforma_tk} and Theorem
\ref{phinormalforma and closeness_tk}, and to show how they are connected to the estimate in Theorem 
\ref{Conjugacy estimate in the TC bifurcation}. Finally we also highlight some issues 
regarding possible generalizations of Theorems 
\ref{Conjugacy estimate in the TC bifurcation} or 
\ref{Conjugacy estimate in the PF bifurcation}
 to higher dimensions.

\begin{figure}[h]
 \centering
 \includegraphics{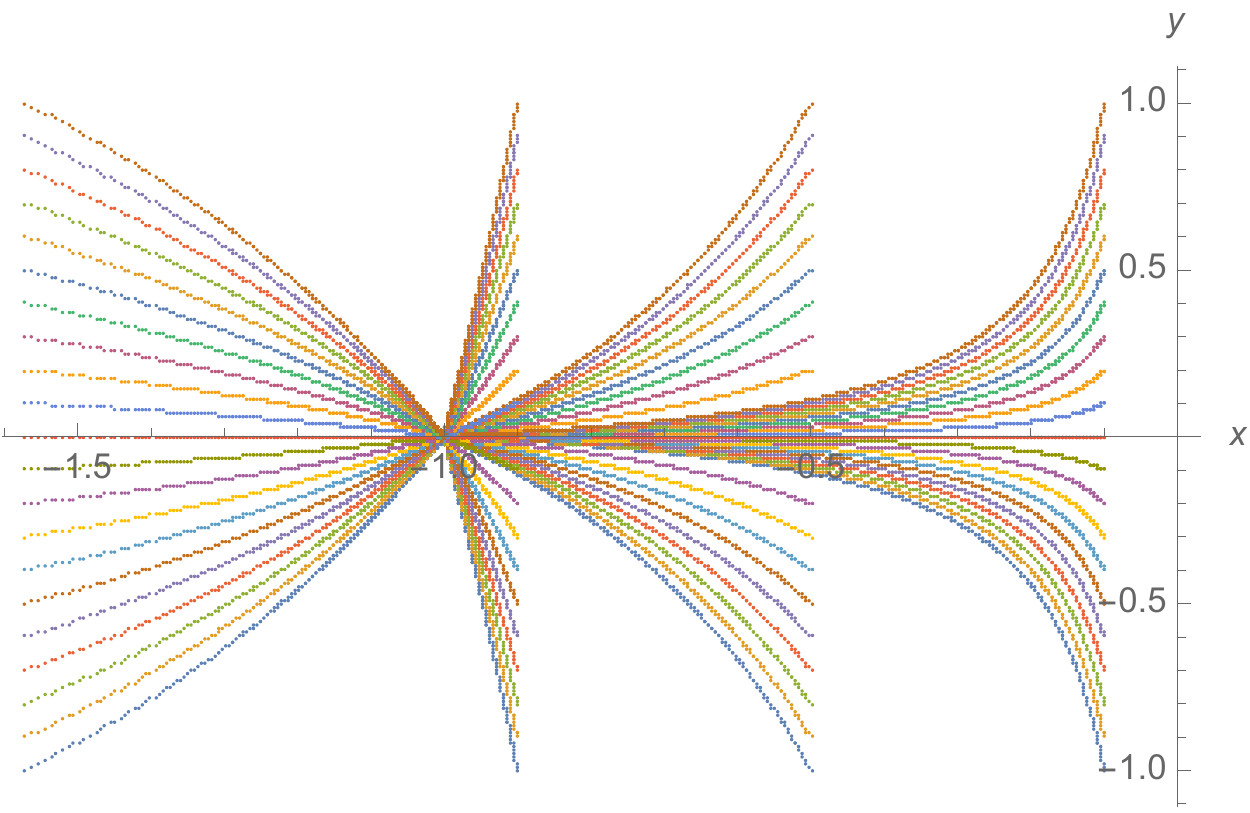}
 \caption{Some iterates of the form
$\left(\Phi\left(n h, x_0,\alpha\right),\Psi\left(n h, y_0,\alpha\right)\right)$ 
corresponding to the time-$h$-map of \eqref{2x2system} with
$h=1/100$ and $\alpha=1$ and started from various initial values $(x_0,y_0)$ show the convergence of the trajectories to the stable fixed point.}\label{figure2Dorbits}
\end{figure}

Let us consider the two-dimensional system
\begin{equation}\label{2x2system}
\begin{cases}
\dot{x}(t)  = &  \alpha x(t)+x^2(t) \\
\dot{y}(t)  = &  -y(t)
\end{cases}
\end{equation}
undergoing a transcritical bifurcation at $\alpha=0$.
By solving \eqref{2x2system},  we see that the time-$h$-map of the exact solution flow $(\Phi,\Psi)$ can be written as   
\begin{equation}\label{timehmapexactsolflow}
\Phi(h,x_0,\alpha) =
\begin{cases}
  \frac{x_0 \alpha e^{\alpha  h}}{\alpha +x_0 \left(1-e^{\alpha  h}\right)} & \text{for } \alpha \neq 0 \\
 \frac{x_0}{1-h x_0} & \text{for } \alpha =0
\end{cases}
\end{equation}
and 
\[
\Psi(h,y_0,\alpha)=y_0 e^{-h}
\]
for any $h\ge 0$, $(x_0,y_0)\in\mathbb{R}^2$ and $\alpha\in\mathbb{R}$ such that the denominators are non-zero. 
Let us apply the classical $4^\mathrm{th}$-order Runge--Kutta method with (some 
small enough) step-size $h>0$  to
\eqref{2x2system}; hence from now on the value of $h$ is fixed. Then the corresponding discretization map $(\varphi,\psi)$ takes the form 
\[\varphi(h,x_0,\alpha)=
\left(1+\alpha  h+\frac{\alpha ^2 h^2}{2}+\frac{\alpha ^3 h^3}{6}+\frac{\alpha ^4 h^4}{24}\right)x_0+\]
\begin{equation}\label{RK4series}
\left(h+\frac{3 \alpha  h^2}{2}+\ldots+ \frac{\alpha ^6 h^7}{96} \right)x_0^2+\ldots +
\frac{h^{15}}{24576}\, x_0^{16} 
\end{equation}
(with intermediate terms suppressed) and
$\psi(h,y_0,\alpha)=\left(1-h+\frac{h^2}{2}-\frac{h^3}{6}+\frac{h^4}{24}\right) y_0$.
By using the semi-group property of the flow, the $n$-fold iterate ($n\in\mathbb{N}$) of $\Phi$ is given by $\left(\Phi(h,\cdot,\alpha)\right)^{[n]}(x_0)$$=\Phi(n h, x_0,\alpha)$, and an analogous
formula holds for $\Psi$. Figure \ref{figure2Dorbits} depicts some iterates 
\[\left(\Phi\left(n\cdot\frac{1}{100}, x_0,1\right),\Psi\left(n\cdot\frac{1}{100}, y_0,1\right)\right)\]
for various initial values $(x_0,y_0)$ and for $0\le n\le 600$. Due to the fact that the system 
\eqref{2x2system} is already decoupled, the center manifolds both for the time-$h$-map
(corresponding to the
fixed point $(0,0)$) and for the discretization map can be
identified with $\mathbb{R}$ for any $\alpha$, see \cite{kuznetsov}. Moreover,
the restriction of the time-$h$-map and the discretization map to this center manifold
can be chosen as $\Phi$ and $\varphi$, respectively.
According to Remark \ref{remark2}, the map $\varphi$ satisfies
the assumptions of Theorem \ref{Conjugacy estimate in the TC bifurcation} with $p=4$.

The qualitative part of Theorem \ref{Conjugacy estimate in the TC bifurcation}---that is, the existence of a conjugacy---implies that the phase portraits induced by the restricted maps
$\Phi$ and $\varphi$ 
near the bifurcation point $(\alpha_0, x_0)=(0,0)$ are topologically the same, and they are both similar to the portrait depicted in Figure \ref{figure3}. 

On the other hand, due to the construction of $J$ presented in Section \ref{tk_constructionoftheconjugacy}, the quantitative part of 
Theorem \ref{Conjugacy estimate in the TC bifurcation}---that is, the estimate \eqref{Theorem2.4closenesspart}---expresses the fact that the orbit of a point $x_0$ under the normal form
${\mathcal{N}_{\Phi}}$ and that of under ${\mathcal{N}_{\vfi}}$ are uniformly close to each other. 
Let us elaborate on this closeness relation. We notice that $\left({\mathcal{N}}_{\Phi}(h,\cdot,\alpha)\right)^{[n]}(x_0)$ and $\left({\mathcal{N}}_{\varphi}(h,\cdot,\alpha)\right)^{[n]}(x_0)$ both converge to the appropriate fixed points as $n\to +\infty$, for any 
$x_0\le 0$ and $\alpha\in \mathbb{R}$ close to $0$, see Figure \ref{figure3}. 
By defining the ``normalized difference''
\begin{equation}\label{deltadef}
\delta(h,x_0,\alpha,n):=\frac{1}{h^p}\left|\left({\mathcal{N}}_{\Phi}(h,\cdot,\alpha)\right)^{[n]}(x_0)-\left({\mathcal{N}}_{\varphi}(h,\cdot,\alpha)\right)^{[n]}(x_0)\right|
\end{equation}
 with $p=4$, the uniform closeness of the orbits means that 
there are constants 
$const>0$, $h_0>0$, $\varepsilon_0>0$ and $\alpha_0>0 $
(depending on the
right-hand side of \eqref{2x2system} and on the chosen discretization method) such that 
\begin{equation}\label{supsupsupsup}
\sup_{h\in(0,h_0]} \ \sup_{x_0\in[-\varepsilon_0, 0]}\  
\sup_{\alpha\in [-\alpha_0,\alpha_0]}\ \sup_{n\in\mathbb{N}} \ \delta(h,x_0,\alpha,n)\le const.
\end{equation}
(For $x_0\in (0,\varepsilon_0]$, one  uses the iterates of the inverses of the normal forms in 
the definition of $\delta$ in \eqref{deltadef}.) 

\begin{figure}[h!]
 \centering
 \includegraphics{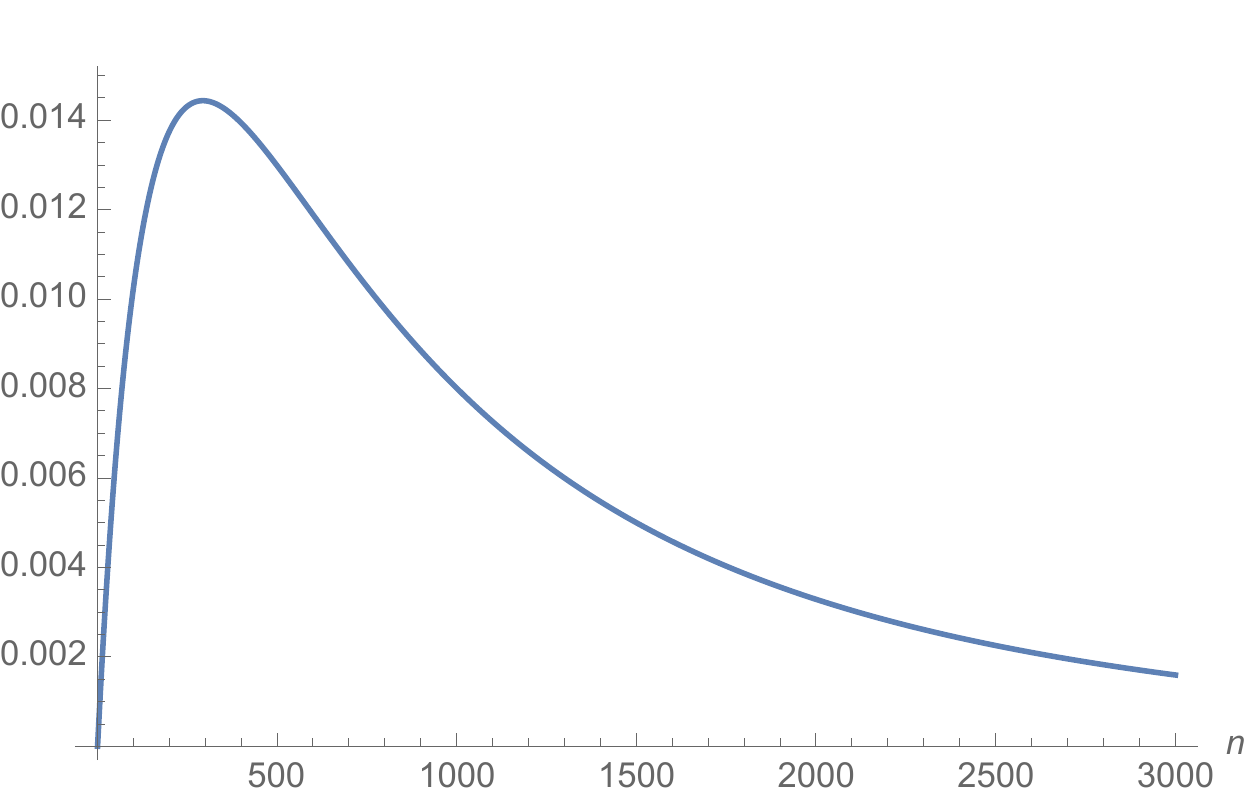}
 \caption{The normalized difference \eqref{normalformsundone} is graphed with $h=1/1000$, $\alpha=-1/2$ and $x_0=-1$.}\label{figurerightorbitdistance}
\end{figure}

\begin{figure}[h!]
 \centering
 \includegraphics{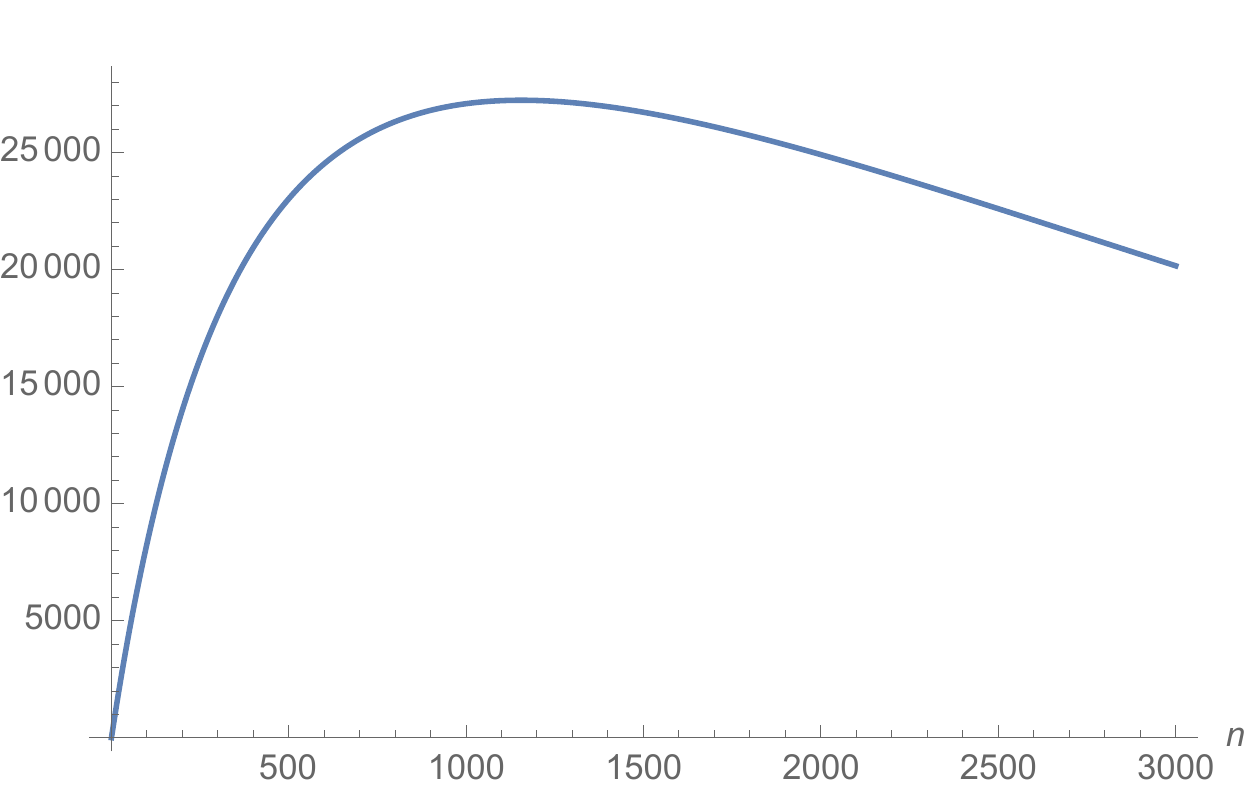}
 \caption{Perturbing the value of $\widetilde\alpha$ in \eqref{normalformsundone}
just by  $10^{-7}$ (and keeping the same $h$, $\alpha$ and $x_0$ values as in Figure \ref{figurerightorbitdistance}) results in a huge increase in the maximum value of \eqref{normalformsundone}.}\label{figureshiftedorbitdistance}
\end{figure}

When casting the above closeness relation in terms of the maps $\Phi$ and $\varphi$ instead of the normal
 forms ${\mathcal{N}}_{\Phi}$ and ${\mathcal{N}}_{\varphi}$, it 
should be emphasized that it is not the orbits
 $\{\left(\Phi(h,\cdot,\alpha)\right)^{[n]}(x_0) : n\in\mathbb{N}\}$ and 
$\{\left(\varphi(h,\cdot,\alpha)\right)^{[n]}(x_0) : n\in\mathbb{N}\}$ that are compared. 
While the step-size $h$ is the same in both cases,  
the initial value $x_0$ and the bifurcation parameter $\alpha$ for (say) the numerical orbit 
generally have to be slightly adjusted to have uniform closeness
between the orbits (see also \cite{beynSIAM}, \cite{garaydisc} or \cite[Section 2.1]{actanumerica}). 
In other words, we prove uniform closeness between the corresponding
members of the sequences
$\{\left(\Phi(h,\cdot,\alpha)\right)^{[n]}(x_0): n\in\mathbb{N}\}$ and 
$\{\left(\varphi(h,\cdot,{\widetilde \alpha})\right)^{[n]}({\widetilde x_0}) : n\in\mathbb{N}\}$, where ${\widetilde x_0}$ and ${\widetilde \alpha}$ are chosen suitably, close to 
$x_0$ and $\alpha$, respectively. (We remark that the choice of 
${\widetilde x_0}$ and ${\widetilde \alpha}$ is in general not unique.)
The normal form transformations in Lemma \ref{Phinormalforma_tk} and Theorem
\ref{phinormalforma and closeness_tk} are preparatory conjugacies,
preceding the construction of the main conjugacy $J$, that provide us with suitable
pairings $x_0 \leftrightarrow {\widetilde x_0}$ and $\alpha \leftrightarrow {\widetilde \alpha}$.
In the present example, these preparatory conjugacies transform 
\eqref{timehmapexactsolflow} into \eqref{tk_normalFi}, and \eqref{RK4series} into 
\eqref{tk_normalfi}. Notice that in general the normal forms have better closeness 
properties 
\eqref{tk_normalformakkozelsege} than the original maps \eqref{kozelseg_tk}.
By considering the proofs of Lemma \ref{Phinormalforma_tk} and Theorem
\ref{phinormalforma and closeness_tk} as given in \cite{dissertation}, one
can explicitly give a formula for ${\widetilde x_0}$ and ${\widetilde \alpha}$:
by inverting the normal form transformations, we get from 
\eqref{supsupsupsup} that
\begin{equation}\label{normalformsundone}
\frac{1}{h^p}\left|\left(\Phi(h,\cdot,\alpha)\right)^{[n]}(x_0)-\frac{1}{\varrho(h,\alpha)}\cdot
\left(\varphi(h,\cdot,\widetilde\alpha(h,\alpha))\right)^{[n]}(\varrho(h,\alpha)\cdot x_0)\right|
\end{equation}
is uniformly bounded in $h\in(0,h_0]$, $x_0\in[-\varepsilon_0, 0]$,
$\alpha\in [-\alpha_0,\alpha_0]$ and $n\in\mathbb{N}$.
The functions $\varrho$ and $\widetilde\alpha$ can explicitly be expressed in terms of
the (inverse functions of some of the) multivariate series expansion coefficient functions of $\Phi$ and $\varphi$, and we have 
$\varrho(h,\alpha)=1+{\cal{O}}(h^p)$ and $\widetilde\alpha(h,\alpha)=\alpha+{\cal{O}}(h^p)$. In 
  our model example
\eqref{timehmapexactsolflow}--\eqref{RK4series}, the functions $\varrho$ and 
$\widetilde\alpha$ contain the exponential function and one root of a quartic polynomial (depending on $h$ and $\alpha$). For simplicity, instead of reproducing 
the exact formulae, here we only give the
first few terms of their series expansion as
\[
\varrho(h,\alpha)=1+\frac{1}{20} \alpha ^4 h^4-\frac{5}{96} \alpha ^5 h^5 +{\cal{O}}\left((\alpha h)^6\right)
\]
and
\[
\widetilde\alpha(h,\alpha)=\alpha\cdot \left(1+\frac{1}{120}\alpha ^4 h^4 -\frac{1}{144} \alpha ^5 h^5+
{\cal{O}}\left((\alpha h)^6\right)\right).
\]
Figure \ref{figurerightorbitdistance} shows the normalized distance \eqref{normalformsundone} with 
$\Phi$ and $\varphi$ given by $\eqref{timehmapexactsolflow}$ and $\eqref{RK4series}$, $h=\frac{1}{1000}$, $\alpha=-\frac{1}{2}$, $x_0=-1$ and $p=4$ for $0\le n \le 3000$. To illustrate the sensitivity of 
\eqref{normalformsundone}  with respect to perturbations, Figure \ref{figureshiftedorbitdistance}
depicts \eqref{normalformsundone} with the same $\Phi$, $\varphi$, $h$, $\alpha$, $x_0$ and $p$, but
instead of $\varphi(h,\cdot,\widetilde\alpha(h,\alpha))$,  the map 
$\varphi(h,\cdot,\widetilde\alpha(h,\alpha)+10^{-7})$ is iterated: 
compare the vertical scale of Figure \ref{figurerightorbitdistance} with that of Figure
\ref{figureshiftedorbitdistance}.

\begin{rem} In \cite{farkas} and \cite{mcli}, conjugacy results
  are proved in the vicinity of a fold bifurcation point in higher dimensions. Both papers show
the existence of a conjugacy between $\Phi(1,\cdot,\alpha)$---i.e., the time-$1$-map of the exact solution operator---and
its numerical approximation \[\left(\varphi(1/N,\cdot,\widetilde\alpha)\right)^{[N]}\]
(with $N\in\mathbb{N}$ sufficiently large). Based on these results it is natural to conjecture
that the existence of a conjugacy between the time-$h$-map of the exact flow and its
numerical approximation can be proved near a TC or a PF bifurcation point 
in higher dimensions as well. 

However, it is not clear \textit{a priori} whether it is possible to prove 
closeness results between the conjugacy and the identity, or what the optimal order ${\cal{O}}(h^q)$  
of a closeness estimate near a TC or PF point would be---one should take into 
account that in general the center manifolds corresponding to 
the continuous and discrete systems do not coincide. We remark that neither \cite{farkas} nor \cite{mcli} contains
any closeness results near a fold bifurcation point. As we pointed out in \cite{garayloczi}, the symmetry argument in  
\cite{farkas} used in the announced ${\cal{O}}(h^p)$ estimate on the center manifold breaks down, hence that
closeness result cannot be considered as proved: as we mentioned in the introduction of 
the present paper, the optimal form of the fold closeness estimate is still an open question
even in one dimension. \quad $\Box$
\end{rem}

\end{document}